\documentclass[11pt]{article}
\def\?{?\vadjust{\vbox to 0pt{\vss\hbox{\kern\hsize\kern1em\large\bf ?!}}}}
\usepackage[T2A]{fontenc}
\usepackage[english]{babel}
\usepackage[tbtags]{amsmath}
\usepackage{amsfonts,amssymb}
\usepackage{mathrsfs}
\usepackage{bm}
\usepackage{esint}
\usepackage{amsthm}

\newcounter{Th}[section]
\newcounter{Lm}[section]
\newcounter{Ca}[section]
\newcounter{ThA}
\newcounter{LmA}
\newcounter{Problem}[section]
\newcounter{Remark}[section]
\newcounter{Example}[section]
\newcounter{Def}[section]
\newcounter{Assum}[section]
\newcounter{Prop}[section]

\def\theTh{\arabic{section}.\arabic{Th}}
\def\theThA{\Alph{ThA}}
\def\theLmA{\Alph{LmA}}
\def\theLm{\arabic{section}.\arabic{Lm}}
\def\theCa{\arabic{section}.\arabic{Ca}}
\def\theProp{\arabic{section}.\arabic{Prop}}

\def\theRemark{\arabic{section}.\arabic{Remark}}
\def\theExample{\arabic{section}.\arabic{Example}}
\def\theDef{\arabic{section}.\arabic{Def}}

\newenvironment{Th}[1][\relax]
    {\medspace\refstepcounter{Th}{\bf Theorem \theTh.}\ \it}
    {\rm\medspace}

\newenvironment{Prop}[1][\relax]
    {\medspace\refstepcounter{Prop}{\bf Proposition \theProp.}\ \it}
    {\rm\medspace}

\newenvironment{ThA}[1][\relax]
    {\medspace\refstepcounter{ThA}{\bf Theorem \theThA.}\ \it}
    {\rm\medspace}

\newenvironment{Lm}[1][\relax]
    {\medspace\refstepcounter{Lm}{\bf Lemma \theLm.}\ \it}
    {\rm\medspace}

\newenvironment{Remark}[1][\relax]
    {\medspace\refstepcounter{Remark}{\bf Remark \theRemark.}\rm\ }
    {\medspace}

\newenvironment{Def}[1][\relax]
    {\medspace\refstepcounter{Def}{\bf Definition \theDef.}\rm\ }
    {\medspace}

\makeatother

\numberwithin{equation}{section}

\date{}

\begin{document}
\author{\bfseries\large A.~I.~Tyulenev
\thanks{Steklov Mathematical Institute of Russian Academy of Sciences (Moscow).
The work of A. I. Tyulenev  was performed at the Steklov International Mathematical Center and supported by the Ministry of Science and Higher Education of the Russian Federation (agreement no. 075-15-2019-1614)
E-mails:
tyulenev-math@yandex.ru, tyulenev@mi.ras.ru}}

\title{Restrictions of Sobolev $W_{p}^{1}(\mathbb{R}^{2})$-spaces to planar rectifiable curves%
}

\maketitle

\begin{abstract}
We construct explicit examples of Frostman-type measures concentrated on arbitrary simple rectifiable curves $\Gamma \subset \mathbb{R}^{2}$ of positive length. Based on such
constructions we obtain for each $p \in (1,\infty)$ an exact description of the trace space $W^{1}_{p}(\mathbb{R}^{2})|_{\Gamma}$ of the first-order Sobolev space $W^{1}_{p}(\mathbb{R}^{2})$ to an arbitrary simple rectifiable curve $\Gamma$
of positive length.
\end{abstract}


\begin{flushleft}
\end{flushleft}
\def\B{\rlap{$\overline B$}B}
{\small\leftskip=10mm\rightskip=\leftskip\noindent}

\section{Introduction}

The problem of the exact description of restrictions of Sobolev spaces $W^{1}_{p}(\mathbb{R}^{n})$, $p \in [1,\infty]$, to different subsets $S \subset \mathbb{R}^{n}$
has a rich history. It takes the origin in the pioneering work of Gagliardo \cite{Gal} where the case  $S=\mathbb{R}^{n-1}$ was considered.
In fact the methods of \cite{Gal} allowed to cover the case when $S$ is a graph of a Lipschitz function $H: \mathbb{R}^{n-1} \to \mathbb{R}$.
Note that this work extended the earlier results by Aronszajn \cite{Ar} and Slobodetskii and Babich \cite{Ba} concerning the
case $p=2$. It should be mentioned that the trace problem for higher order Sobolev spaces $W^{m}_{p}(\mathbb{R}^{n})$, $p \in (1,\infty)$, $m \in \mathbb{N}$
in the case $S=\mathbb{R}^{d}$, $d \in [1,n-1] \cap \mathbb{N}$ was covered by Besov in the fundamental paper \cite{Besov}.

\textit{In the case $p=\infty$}, the Sobolev space $W^{1}_{\infty}(\mathbb{R}^{n})$ can be identified with the space
$\operatorname{LIP}(\mathbb{R}^{n})$  of Lipschitz functions on $\mathbb{R}^{n}$ and it is known (see McShane-Whitney extension lemma in section 4.1 of \cite{Kos}) that for
any closed set $S \subset \mathbb{R}^{n}$ the restriction $\operatorname{LIP}(\mathbb{R}^{n})|_{S}$
coincides with the space $\operatorname{LIP}(S)$ of Lipschitz functions on $S$ and that, furthermore, the
classical Whitney extension operator linearly and continuously maps the space $\operatorname{LIP}(S)$
into the space  $\operatorname{LIP}(\mathbb{R}^{n})$  (see e.g., \cite{St}, Chapter 6).

\textit{In the case} $p=1$ much less is known. Indeed, as it was shown in \cite{Gal} for
the case $S=\mathbb{R}^{n-1}$ the trace space on $\mathbb{R}^{n-1}$ of the Sobolev space $W_{1}^{1}(\mathbb{R}^{n})$ can be identified
with $L_{1}(\mathbb{R}^{n-1})$ as a linear space, the corresponding norms being equivalent. However,
the extension operator constructed by Gagliardo is nonlinear. Furthermore, it was shown by J. Peetre \cite{Pet}
(see also section 5 in \cite{PelWoj}) that  any bounded map from $L_{1}(\mathbb{R}^{n-1})$ to $W_{1}^{1}(\mathbb{R}^{n})$
which is right inverse to the trace map is nonlinear.

In the sequel we deal with the case $p \in (1,\infty)$ only.

After \cite{Gal} a big progress was made by several mathematicians \cite{Jon2}, \cite{Shv1},  \cite{Shv2}, \cite{Ihn}, \cite{TrW}, \cite{Kal},
\cite{Ry}
in the direction of relaxation of extra assumptions on the sets $S$.

We have to note that in \cite{Shv2} the corresponding trace problem was solved without any assumptions on $S$.
However, in that paper only the case $p > n$ was considered. This case is special and exploits techniques different from that
of used in other papers mentioned above. Unfortunately, such techniques do not allow to attack the case $p \in (1,n]$.

Recall that given $d \in (0,n]$, a closed set $S \subset \mathbb{R}^{n}$ is said to be \textit{Ahlfors-David $d$-regular} provided that there exist
constants $c^{1}_{S},c^{2}_{S} > 0$ such that
\begin{equation}
\label{eq1.1}
c^{1}_{S}r^{d} \le  \mathcal{H}^{d}(Q(x,r) \cap S) \le c^{2}_{S}r^{d} \quad \hbox{for every} \quad x \in S, \quad r \in (0,1].
\end{equation}
In \eqref{eq1.1} we set $Q(x,r):=\prod_{i=1}^{n}[x_{i}-r,x_{i}+r]$ and by $\mathcal{H}^{d}$ we denoted the $d$-Hausdorff measure.
We will also call condition \eqref{eq1.1} the Ahlfors-David $d$-regularity condition.

Summarizing results and methods of papers  \cite{Jon2}, \cite{Shv1}, \cite{Ihn} restricted to the case of the
first-order Sobolev spaces $W_{p}^{1}(\mathbb{R}^{n})$ one can obtain for any fixed $d \in (0,n]$ and $p \in (\max\{1,n-d\},\infty)$
an exact description of the trace space of the space $W_{p}^{1}(\mathbb{R}^{n})$ to any closed set $S \subset \mathbb{R}^{n}$ satisfying the Ahlfors-David
$d$-regularity condition \eqref{eq1.1}.

Recently Rychkov introduced \cite{Ry} the concept of $d$-thick sets. Recall that given $d \in [0,n]$, a set $S \subset \mathbb{R}^{n}$
is said to be \textit{$d$-thick} if there exists a constant $c^{3}_{S} > 0$ such that
\begin{equation}
\label{eq1.2}
c^{3}_{S}r^{d} \le  \mathcal{H}^{d}_{\infty}(Q(x,r) \cap S) \quad \hbox{for every} \quad x \in S, \quad r \in (0,1],
\end{equation}
where by $\mathcal{H}^{d}_{\infty}$ we denoted the so-called $d$-Hausdorff content.
We have to note that condition \eqref{eq1.2} is much weaker than \eqref{eq1.1}. It was noted in \cite{Ry} and proved in \cite{TV}
that every path-connected set $S \subset \mathbb{R}^{n}$ consisting of more than one point is $1$-thick.
It is clear that a generic path-connected set $S$ fails to satisfy Ahlfors-David $1$-regularity
condition. In \cite{Ry} trace criteria for Besov $B^{s}_{p,q}(\mathbb{R}^{n})$ and Lizorkin-Triebel $F^{s}_{p,q}(\mathbb{R}^{n})$ spaces were
obtained for $d$-thick sets $S$. However, that criteria
were not fully intrinsic and were based on atomic-type characterizations of function spaces.
Furthermore, in the case $s \in \mathbb{Z}$ the extra assumption $d > n-1$ was required. 
In particular, that restrictions did not allow to attack the trace problem for Sobolev spaces $W^{1}_{p}(\mathbb{R}^{n})$ in the case
of $d$-thick sets $S$ with $d \in [0,n-1]$. 

We have to mention papers \cite{TrW}, \cite{Kal}. In \cite{Kal} the trace problem for Sobolev spaces was considered
in the case when $S \subset \mathbb{R}^{2}$ is a single cusp satisfying some extra regularity assumptions. In \cite{TrW} the trace problem for Besov $B^{s}_{p,q}(\mathbb{R}^{n})$ and Lizorkin-Triebel $F^{s}_{p,q}(\mathbb{R}^{n})$ spaces was considered in the case when $S \subset \mathbb{R}^{n}$ is a domain satisfying the so-called internal and external regularity assumptions. However, the corresponding criteria were not fully intrinsic and used atomic-type characterizations.

Very recently \cite{TV}, given parameters $d \in [0,n]$, $p \in (\max\{1,n-d\},\infty)$ and a closed $d$-thick set $S \subset \mathbb{R}^{n}$ an exact description of traces of functions $F \in W^{1}_{p}(\mathbb{R}^{n})$ to the set $S$ was obtained. As far as we know it was a first result concerning trace problems
for the first-order Sobolev spaces $W^{1}_{p}(\mathbb{R}^{n})$ obtained for the range $p \in (1,n]$ in such a high generality.
Furthermore, it is possible to generalize that results to the case of weighted Sobolev spaces \cite{TV2}.

Analysis of results obtained in \cite{Jon2}, \cite{Shv1}, \cite{Ihn} shows that in the case of Ahlofrs-David $d$-regular sets $S$ the
only measure which played a crucial role in the solution of trace problems is the $d$-Hausdorff measure $\mathcal{H}^{d}$.
But this is not the case for $d$-thick sets.
Unfortunately, in this case  one has to work with a sequence of measures with prescribed growth conditions instead of the only "nice measure".
Recall \cite{TV} that given a closed $d$-thick set $S \subset \mathbb{R}^{n}$, there exists a sequence of Radon measures $\{\mathfrak{m}_{k}\}:=\{\mathfrak{m}_{k}\}_{k \in \mathbb{N}_{0}}$ with the following properties:


{\rm (1)} for every $k \in \mathbb{N}_{0}$
\begin{equation}
\operatorname{supp}\mathfrak{m}_{k}=S;
\end{equation}

{\rm (2)} there exists a constant $C^{1} > 0$ such that for each $k \in \mathbb{N}_{0}$
\begin{equation}
\label{eq2.15}
\mathfrak{m}_{k}(Q(x,r)) \le C^{1} r^{d} \quad  \text{ for every } x \in \mathbb{R}^{n}
\text{ and every } r \in (0,2^{-k}];
\end{equation}

{\rm (3)} there exists a constant $C^{2} > 0$ such that for each $k \in \mathbb{N}_{0}$
\begin{equation}
\label{eq2.16}
\mathfrak{m}_{k}(Q(x,2^{-k})) \geq C^{2} 2^{-dk} \quad \text{for every} \quad x \in S;
\end{equation}

{\rm (4)} for each $k \in \mathbb{N}_{0}$ the measure $\mathfrak{m}_{k}=w_{k}\mathfrak{m}_{0}$ with $w_{k} \in L_{\infty}(\mathfrak{m}_{0})$ and
\begin{equation}
\label{eq2.17}
2^{d-n}w_{k+1}(x) \le w_{k}(x) \le w_{k+1}(x) \quad  \text{for} \quad \mathfrak{m}_{0}-a.e. \quad x \in S.
\end{equation}

In what follows we call any sequence of Radon measures $\{\mathfrak{m}_{k}\}$ satisfying items (1)--(4) above
\textit{a $d$-regular on $S$.}

We denote by $C^{1}_{\{\mathfrak{m}_{k}\}}$ and $C^{2}_{\{\mathfrak{m}_{k}\}}$ the minimum among constants $C^{1}$ in \eqref{eq2.15} and
the maximum among constants $C^{2}$  in \eqref{eq2.16}  respectively.

In order to describe the aim of this paper and for the sake of completeness of exposition
we have to formulate one particular case of the main result obtained in \cite{TV}.
In doing so we firstly recall some basic concepts which were introduced in \cite{TV}.
Given an arbitrary sequence $\{\mathfrak{m}_{k}\}$ of Radon measures on $\mathbb{R}^{n}$, we define for each $t \in (0,1]$
the \textit{Calderon-type maximal function} of $f$ with respect to $\{\mathfrak{m}_{k}\}$ as
\begin{equation}
\label{eq1.3}
f^{\sharp}_{\{\mathfrak{m}_{k}\}}(x,t):=\sup\limits_{k \in \mathbb{N}_{0}, 2^{-k} \geq t}
2^{k}\inf\limits_{c \in \mathbb{R}}\fint\limits_{Q(x,2^{-k})}|f(y)-c|\,d\mathfrak{m}_{k}(y),
\end{equation}
where the corresponding averaged integrals are assumed to be zero in the case $\mathfrak{m}_{k}(Q(x,2^{-k}))=0$.

Finally, given a closed set $S \subset \mathbb{R}^{n}$ and a parameter $\lambda \in (0,1)$ we define for each $t \in (0,1]$
the \textit{maximal $\lambda$-porous at the scale $t$ subset of $S$} as
\begin{equation}
\label{eq1.7'}
S_{t}(\lambda):=\{x \in S: \hbox{there exists } y \in Q(x,t) \hbox{ s.t. } Q(y,\lambda t) \subset \mathbb{R}^{n} \setminus S\}.
\end{equation}

Given  parameters $d \in (0,n]$, $p \in (\max\{1,n-d\},\infty)$ and a closed set $S$ with $\mathcal{H}^{d}(S) > 0$,
by the symbol $W_{p}^{1}(\mathbb{R}^{n})|_{S}$ we denote the trace space
of the Sobolev space $W^{1}_{p}(\mathbb{R}^{n})$ to the set $S$ (see the next section for the precise definition).

Recall briefly the construction of the extension operator from \cite{TV}.
Given $d \in (0,n]$ and a  $d$-thick closed set  $S \subset \mathbb{R}^{n}$, let $\{Q_{\alpha}\}_{\alpha \in I}$ be the Whitney
decomposition of $\mathbb{R}^{n} \setminus S$ and let $\mathcal{I} \subset I$ be the index set corresponding
to the Whitney cubes with side lengthes $\le 1$. Let $\{\varphi_{\alpha}\}_{\alpha \in I}$ be the corresponding partition of unity (see \cite{TV} for details).
For any cube $Q_{\alpha}=Q(x_{\alpha},r_{\alpha})$, $\alpha \in I$ we define the cube
$\widetilde{Q}_{\alpha}:=Q(\widetilde{x}_{\alpha},r_{\alpha})$, where $\widetilde{x}_{\alpha}$ is an arbitrary metric projection
of $x_{\alpha}$ to the set $S$. Let $\{\mathfrak{m}_{k}\}$ be an arbitrary $d$-regular on $S$ sequence of measures.
Let $f \in L^{\operatorname{loc}}_{1}(\mathfrak{m}_{k})$ for some (and hence every)
$k \in \mathbb{N}_{0}$. We set $k(r):=[\log_{2}r^{-1}]$ and define
\begin{equation}
\label{eq2.37'}
F(x)=\operatorname{Ext}_{S,\{\mathfrak{m}_{k}\}}[f](x):=\sum\limits_{\alpha \in \mathcal{I}}\varphi_{\alpha}(x)\fint\limits_{\widetilde{Q}_{\alpha} \cap S}f(\widetilde{x})
\,d\mathfrak{m}_{k(r_{\alpha})}(\widetilde{x})
,\quad x \in \mathbb{R}^{n}.
\end{equation}

\begin{ThA}
\label{ThA}
Let $d \in (0,n)$ and $p \in (\max\{1,n-d\},\infty)$. Let $S \subset \mathbb{R}^{n}$ be a closed $d$-thick set with $\mathcal{H}^{n}(S)=0$.
Given a  $d$-regular sequence of measures $\{\mathfrak{m}_{k}\}:=\{\mathfrak{m}_{k}\}_{k \in \mathbb{N}_{0}}$  on $S$, a measurable function $f: S \to \mathbb{R}$ belongs to the trace space $W_{p}^{1}(\mathbb{R}^{n})|_{S}$ if and only if
the following conditions hold:

{\rm (1)} for $\mathcal{H}^{d}$-almost every point $x \in S$ it holds
\begin{equation}
\label{eq1.3'}
\lim\limits_{k \to \infty}\fint\limits_{Q(x,2^{-k}) \cap S} |f(x)-f(z)|\,d\mathfrak{m}_{k}(z) = 0;
\end{equation}

{\rm (2)} there exists a number $\lambda_{0} \in (0,1)$ such that (we set $S_{k}(\lambda):=S_{2^{-k}}(\lambda)$)
\begin{equation}
\label{eq1.9}
\mathcal{BN}_{\{\mathfrak{m}_{k}\},p,\lambda_{0}}[f]:=\|f|L_{p}(\mathfrak{m}_{0})\|+ \Bigl(\sum_{k=1}^{\infty}2^{k(d-n)}\int\limits_{S_{k}(\lambda_{0})}\bigl(f^{\sharp}_{\{\mathfrak{m}_{k}\}}(x,2^{-k})\bigr)^{p}\,d\mathfrak{m}_{k}(x)\Bigr)^{\frac{1}{p}} < +\infty.
\end{equation}

Furthermore, for every $\lambda \in (0,\lambda_{0}]$
there exists a constant $C > 0$ depending only on $p,n,d,C^{1}_{\{\mathfrak{m}_{k}\}},C^{2}_{\{\mathfrak{m}_{k}\}}$ and $\lambda$ such that
\begin{equation}
\label{eq1.4}
\frac{1}{C}\mathcal{BN}_{\{\mathfrak{m}_{k}\},p,\lambda}[f] \le \|f|W_{p}^{1}(\mathbb{R}^{n})|_{S}\|  \le  C \mathcal{BN}_{\{\mathfrak{m}_{k}\},p,\lambda}[f].
\end{equation}
The extension operator $\operatorname{Ext}_{S,\{\mathfrak{m}_{k}\}}$ is a right inverse operator
for the usual trace operator. It  maps the trace
space $W_{p}^{1}(\mathbb{R}^{n})|_{S}$ to the space $W_{p}^{1}(\mathbb{R}^{n})$ linearly and continuously.
\end{ThA}

In fact Theorem \ref{ThA} is an almost immediate consequence of Theorem 2.1 from \cite{TV}. The only delicate point
which we have to mention is that condition \eqref{eq1.3'} is more rough than the corresponding condition
in \cite{TV} because it appeals to Hausdorff measures instead of capacities.
On the other hand, our definition of the trace space given in this paper also appeals to Hausdorff measures rather than capacities (compare Definition \ref{Def2.1} in this paper
with Definitions 2.7, 2.8, 2.9 in  \cite{TV}).
Hence, a careful analysis of proofs of Lemma 4.3 and Theorem 4.2 in \cite{TV} together with very well known
relations between capacities and Hausdorff measures shows that the using of \eqref{eq1.3'} is justified.

In practice the criterion given in Theorem \ref{ThA} is not so easy to verify by the following reasons:

{\rm (1)} it is difficult to check a delicate condition \eqref{eq1.3'};

{\rm (2)} in fact the typical construction of a $d$-regular sequence of measures $\{\mathfrak{m}_{k}\}$ given in \cite{TV}
was not fully constructive. Indeed, we used the classical Frostman-type arguments which are based on an inductive algorithm and on the weak limit
procedure. In practice it is not so easy to present explicit expressions for such measures.
Hence, the construction of explicit examples of
$d$-regular sequences of measures $\{\mathfrak{m}_{k}\}$ on different closed $d$-thick sets $S \subset \mathbb{R}^{n}$ is of great importance.

The aim of this paper is to demonstrate that one can overcome the difficulties described in items (1) and (2) above
and simplify Theorem A in the case when $S$ is \textit{a planar simple (i.e. without self-intersections) rectifiable curve
$\Gamma \subset \mathbb{R}^{2}$ with positive length}. In this case $\Gamma$ is a $1$-thick set. We construct \textit{the special $1$-regular on $\Gamma$
sequence of measures} $\{\mu_{k}[\Gamma]\}_{k \in \mathbb{N}_{0}}$. The advantage of the sequence $\{\mu_{k}[\Gamma]\}_{k \in \mathbb{N}_{0}}$ is that the
measures $\mu_{k}[\Gamma]$, $k \in \mathbb{N}_{0}$ have explicit expressions.
Furthermore, having at disposal such measures one can get rid of a sophisticated condition \eqref{eq1.3'} in the statement of Theorem
\ref{ThA}.
We have to note that results obtained in this paper \textit{are new and could not be obtained by previously known methods.}
Indeed, as far as we know \textit{explicit constructions of Frostman-type measures} on arbitrary planar rectifiable curves $\Gamma$ were not
considered in the literature before. Furthermore, even in the the case of a planar simple rectifiable curve
$S=\Gamma \subset \mathbb{R}^{2}$ of positive length results formulated in Theorem A could not be obtained by previously known methods.
Indeed, it was unaccessible by methods of \cite{Ry} because
in the case of Sobolev spaces the corresponding trace problem was considered only for $d$-thick sets $S \subset \mathbb{R}^{n}$ satisfying the additional
requirement $d > n-1$. On the other hand, there are planar rectifiable curves of positive length that fail to satisfy the Ahlfors-David $1$-regularity
condition (and hence, fail to satisfy the Ahlfors-David $d$-regularity conditions for all $d \in (0,2]$) and hence do not fall into the scope of \cite{Ihn}.

Let us informally explain why in the present paper we restrict ourselves to simple rectifiable curves $\Gamma \subset \mathbb{R}^{2}$
only.

{\rm (A)} In fact arguments of section 3 work for simple rectifiable curves in $\mathbb{R}^{n}$, $n \geq 2$
with positive length. The only minor problem is that in the case $n > 2$ the corresponding expressions for measures will be more technical.
On the other hand, what is more important is the application of Theorem 3.1 in section 4 where
restriction to the dimension $n \in \mathbb{N}$ of the ambient space becomes essential. More precisely, as we have already mentioned any rectifiable curve of positive length
is a $1$-thick subset of $\mathbb{R}^{n}$. But Theorem 4.1 (which is a keystone for the main result) works for $d$-thick sets
with $d \in [n-1,n]$. This obstruction justifies the working in the 2-dimensional plane.

{\rm (B)} The main reason why we restrict ourselves to the case of curves instead of general $1$-thick sets in $\mathbb{R}^{2}$ is
that the corresponding expressions for the Frostman-type measures concentrated on general 1-thick sets will be much less transparent. Roughly speaking the main technical
advantage of our construction which works for curves is the reduction of $1$-dimensional Frostman-type measures to the $0$-dimensional Frostman-type 
measures.
Indeed, it is well known that any simple planar rectifiable curve $\Gamma$
have finite number of intersections with "almost every" line parallel to the coordinate axes. This allows
to built the corresponding $0$-dimensional Frostman-type measure with a help of elementary combinatorial
arguments. Clearly, one can not hope to make a similar trick for general $1$-thick sets in $\mathbb{R}^{2}$ because intersections
of such sets with lines can have a complicated geometry.
Indeed, the general $1$-thick sets can be composed of pieces with different dimensions.

\textbf{Acknowledgements.}
I dedicate this paper to my "scientific grandfather" academician of Russian Academy of Science, Professor S. M. Nikol'skii. The paper was written in 2020,
the year of his 115-anniversary.

I thank P. Shvartsman who read my joint paper with S. Vodop'yanov \cite{TV}
and inspired me to built concrete examples of Frostman-type measures.

I am grateful to V. Bogachev and N. Gigli for the fruitful discussions which helped to clarify
some ideas of this paper.
Furthermore, I am grateful to  D. Stolyarov and A. Volberg for the valuable remarks.

I would like to thank the referees for very careful reading and numerous comments and suggestions, which led to
improvements of the manuscript.

\section{Preliminaries}

Throughout the paper $C,C_{1},C_{2},...$ will be generic positive
constants. These constants can change even
in a single string of estimates. The dependence of a constant on certain parameters is expressed, for example, by
the notation $C=C(n,p,k)$. We write $A \approx B$ if there is a constant $C \geq 1$ such that $A/C \le B \le C A$.
Given a number $c \in \mathbb{R}$ we denote by $[c]$ the integer part of $c$.

\subsection{Geometric Measure Theory Background}
If no otherwise stated we let $\mathbb{R}^{n}$, $n \geq 1$ denote the linear space of all strings $x=(x_{1},...,x_{n})$ of real numbers \textit{equipped  with the uniform  norm} $\|\cdot\|_{\infty}$, i.e. $\|x\|_{\infty}:=\max\{|x_{1}|,...,|x_{n}|\}$. As usual $\overline{\mathbb{R}}:=\mathbb{R} \cup \{-\infty\} \cup \{+\infty\}$. Given a set $E \subset \mathbb{R}^{n}$ we denote by $\operatorname{int}E$, $\overline{E}$ and $E^{c}$
the interior, the closure and the complement (in $\mathbb{R}^{n}$) of $E$ respectively. Given a set $E \subset \mathbb{R}^{n}$ we will always denote by
$\chi_{E}$ the characteristic function of $E$. By \textit{a cube} $Q$ in $\mathbb{R}^{n}$ we mean a \textit{closed} cube with sides parallel to the coordinate axes.
We say that $E \subset \mathbb{R}^{n}$ is  \textit{a measurable set} if $E$ to the standard $\sigma$-algebra off all Lebesgue measurable sets.
Given a measurable set $E \subset \mathbb{R}^{n}$, we say that a function $f: E \to \mathbb{R}$ is \textit{measurable} if for any $c \in \mathbb{R}$ the set $f^{-1}((c,+\infty])$
is measurable.

In the sequel given a metric space $\operatorname{X}=(\operatorname{X},\operatorname{d})$, by a measure on $\operatorname{X}$ we mean only \textit{a nonnegative Borel
measure} on $\operatorname{X}$. Given a measure $\mathfrak{m}$ on $\operatorname{X}$ and a nonempty Borel set $S \subset \operatorname{X}$, we define the \textit{restriction $\mathfrak{m}\lfloor_{S}$ of} $\mathfrak{m}$ to $S$ as
\begin{equation}
\mathfrak{m}\lfloor_S (E) := \mathfrak{m}(E \cap S) \quad \hbox{for any Borel set} \quad E \subset \operatorname{X}.
\end{equation}
Given two metric spaces $(\operatorname{X}_{1},\operatorname{d}_{1})$, $(\operatorname{X}_{2},\operatorname{d}_{2})$, Borel map $G: \operatorname{X}_{1} \to \operatorname{X}_{2}$
and Borel measure $\mathfrak{m}$ on $\operatorname{X}_{1}$ we define the \textit{push-froward measure} $G_{\sharp}\mathfrak{m}$ on $\operatorname{X}_{2}$ by the equality
\begin{equation}
G_{\sharp}\mathfrak{m}(E):=\mathfrak{m}(G^{-1}(E)) \quad \hbox{for any Borel set} \quad E \subset \operatorname{X}_{2}.
\end{equation}
Let $\mathfrak{m}$ be an arbitrary measure on $\mathbb{R}^{n}$. Given $f \in L^{\text{\rm loc}}_{1}(\mathbb{R}^{2},\mathfrak{m})$, we set for every
Borel set $G \subset \mathbb{R}^{n}$ with $\mathfrak{m}(G) < +\infty$
\begin{equation}
\label{eq.average}
\fint\limits_{G}f(x)\,d\mathfrak{m}(x):=
\begin{cases}
&\frac{1}{\mathfrak{m}(G)}\int\limits_{G}f(x)\,d\mathfrak{m}(x), \quad \hbox{if} \quad \mathfrak{m}(G) > 0;\\
&0, \quad  \hbox{if} \quad \mathfrak{m}(G) = 0.
\end{cases}
\end{equation}
Given a Radon measure $\mathfrak{m}$ on $\mathbb{R}^{n}$, we set for every cube $Q \subset \mathbb{R}^{n}$
\begin{equation}
\widetilde{\mathcal{E}}_{\mathfrak{m}}[f](Q):=\fint\limits_{Q}\fint\limits_{Q}|f(y)-f(z)|\,d\mathfrak{m}(y)d\mathfrak{m}(z).
\end{equation}

Recall that Calderon-type maximal functions $f^{\sharp}_{\{\mathfrak{m}_{k}\}}$ were defined in the introduction.

\begin{Prop}
\label{Prop23}
Let $\{\mathfrak{m}_{k}\}=\{\mathfrak{m}_{k}\}_{k \in \mathbb{N}_{0}}$ be a sequence of Radon measures on $\mathbb{R}^{n}$.
Then for each $t \in (0,1]$
and every $x \in \mathbb{R}^{n}$ it holds
\begin{equation}
\begin{split}
\label{eq231}
&\frac{1}{2}\sup\limits_{k \in \mathbb{N}_{0}, 2^{-k} \geq t}2^{k}\widetilde{\mathcal{E}}_{\mathfrak{m}_{k}}[f](Q(x,2^{-k})) \le f^{\sharp}_{\{\mathfrak{m}_{k}\}}(x,t) \le
\sup\limits_{k \in \mathbb{N}_{0}, 2^{-k} \geq t}2^{k}\widetilde{\mathcal{E}}_{\mathfrak{m}_{k}}[f](Q(x,2^{-k})).
\end{split}
\end{equation}
\end{Prop}

\begin{proof}
We fix arbitrary $x \in \mathbb{R}^{n}$,$r \in (0,1]$ and set $Q=Q(x,r)$ for brevity. The first inequality in \eqref{eq231} follows form the fact that for each $k \in \mathbb{N}_{0}$ and any constant $c \in \mathbb{R}$ we have
\begin{equation}
\begin{split}
\notag
&\widetilde{\mathcal{E}}_{\mathfrak{m}_{k}}[f](Q) \le \fint\limits_{Q}|f(y)-c|\,d\mathfrak{m}_{k}(y)
+\fint\limits_{Q}|-f(z)+c|\,d\mathfrak{m}_{k}(z)\\
& \le 2\fint\limits_{Q}|f(y)-c|\,d\mathfrak{m}_{k}(y).
\end{split}
\end{equation}
The second inequality in  \eqref{eq231} follows from the estimate
\begin{equation}
\begin{split}
\notag
&\inf\limits_{c \in \mathbb{R}}\fint\limits_{Q}|f(y)-c|\,d\mathfrak{m}_{k}(y)
\le \fint\limits_{Q}\Bigl|f(y)-\fint\limits_{Q}f(z)\,d\mathfrak{m}_{k}(z)\Bigr|\,d\mathfrak{m}_{k}(y) \le \widetilde{\mathcal{E}}_{\mathfrak{m}_{k}}[f](Q), \quad k \in \mathbb{N}_{0}.
\end{split}
\end{equation}
\end{proof}

In this paper it will be convenient (in the sense that there are several definitions of Hausdorff measures
in the literature that give the same values up to some universal constants) to follow \cite{Falc} (see section 1.2 therein) and define for each $d \in [0,n]$ and $\delta \in (0,\infty]$
\begin{equation}
\notag
\mathcal{H}^{d}_{\delta}(E):=\inf\sum\limits_{i=1}^{\infty}(\operatorname{diam}U_{i})^{d},
\end{equation}
where the infimum is over all coverings $\{U_{i}\}_{i \in \mathbb{N}}$ of $E$ with $\operatorname{diam}U_{i} < \delta$, $i \in \mathbb{N}$.
Now,  the \textit{$d$-Hausdorff measure $\mathcal{H}^{d}$ of the set} $E$
is defined as $\mathcal{H}^{d}(E):=\lim_{\delta \to 0}\mathcal{H}^{d}_{\delta}(E)$.
By the \textit{$d$-Hausdorff content of $E$} we mean $\mathcal{H}^{d}_{\infty}(E)$.

The following proposition is an immediate consequence of \eqref{eq2.15} and definition of the
measure $\mathcal{H}^{d}$. We omit an elementary proof.

\begin{Prop}
\label{Prop2.5}
Let $S \subset \mathbb{R}^{n}$ be a closed $d$-thick set for some $d \in (0,n]$. Let $\{\mathfrak{m}_{k}\}$
be a $d$-regular sequence of measures on $S$. Then, for every $k \in \mathbb{N}_{0}$ the measure $\mathfrak{m}_{k}$
is absolutely continuous with respect to $\mathcal{H}^{d}\lfloor_{S}$. Furthermore, for each $k \in \mathbb{N}_{0}$ it holds
\begin{equation}
\label{eq26'''}
\mathfrak{m}_{k}(E) \le  C^{1}_{\{\mathfrak{m}_{k}\}} \mathcal{H}^{d}(E) \quad \hbox{for any Borel set} \quad E \subset S.
\end{equation}
\end{Prop}

\begin{Remark}
Note that the right hand side of \eqref{eq26'''} can be equal $+\infty$ and hence the corresponding estimate
is trivial.
\end{Remark}

Given a measurable set $E \subset \mathbb{R}^{n}$, recall that a map $f: E \to \mathbb{R}$, $E \subset \mathbb{R}$ is said to have \textit{the Lusin property}
if  for any set $E_{0} \subset E$ of Lebesgue measure zero the image $f(E_{0})$ has Lebesgue measure zero.
The following result is a particular case of Theorem 12 of \cite{Haj} (see also Theorem 4.3.3 in \cite{Lin} where minor modifications are required). In fact it will be a keystone
for the whole section 3 below.

\begin{Prop}
\label{Prop2.2}
Let a map $\Phi: \mathbb{R} \to \mathbb{R}$ be continuous, differentiable almost everywhere and has Lusin property.
Let $g: \mathbb{R} \to [0,+\infty]$ be a measurable function.
Then for any measurable set $E\subset\mathbb{R}$
\begin{equation}
\label{eq2.8}
\int\limits_{E}g(x)|\Phi'(x)|\,dx = \int\limits_{\Phi(E)} \sum\limits_{x \in \Phi^{-1}(y) \cap E}g(x)\,d\mathcal{H}^{1}(y)=
\int\limits_{\Phi(E)}\Bigl(\int\limits_{\Phi^{-1}(y)\cap E}g(x)\,d\mathcal{H}^{0}(x)\Bigr)\,d\mathcal{H}^{1}(y).
\end{equation}
\end{Prop}

Let $S \subset \mathbb{R}^{n}$ be a closed set. We recall that given $\lambda \in (0,1)$
and $t \in (0,1]$, the definition of the set $S_{t}(\lambda)$ was given in the introduction.
In what follows for every $k \in \mathbb{N}$ we set $S_{k}(\lambda):=S_{2^{-k}}(\lambda)$.
It is natural to ask whether there exists a parameter $\lambda \in (0,1)$ such that the union $\cup_{k \in \mathbb{N}_{0}}S_{k}(\lambda)$ contains $S$ or equivalently
$\max_{k \in \mathbb{N}_{0}}\chi_{S_{k}(\lambda)}(x)=1$ for every $x \in S$?
Unfortunately, this is not the case in general.
In fact in the sequel sets with even a more restrictive properties will be important for us.
More precisely, we introduce the following concept.

\begin{Def}
\label{Def.quasipor}
Let $d \in (0,n]$ and $\lambda \in (0,1)$. Let $S \subset \mathbb{R}^{n}$ be a closed set with $\mathcal{H}^{d}(S) > 0$. We say that $S$
is \textit{$(d,\lambda)$-quasi-porous} if
\begin{equation}
\lim_{k \to \infty}\chi_{S_{k}(\lambda)}(x)=1 \quad \hbox{for} \quad \mathcal{H}^{d}-\hbox{a.e.} \quad x \in S.
\end{equation}
\end{Def}

\begin{Remark}
\label{Rem.porous}
It is obvious that if a set $S$ is $(d,\lambda)$-quasi-porous for some $\lambda \in (0,1)$, then it is $(d,\lambda')$-quasi-porous
for every $\lambda' \in (0,\lambda]$.
\end{Remark}

The following lemma gives one useful and simple sufficient conditions for a given closed set $S$ to be $(d,\lambda)$-quasi-porous.

\begin{Lm}
\label{Lm.porous}
Let $d \in (0,n)$. Suppose that a closed set $S \subset \mathbb{R}^{n}$ is such that:

{\rm (1)} $\mathcal{H}^{d}(S) \in (0,\infty)$;

{\rm (2)} there is a constant $c_{1} > 0$  such that
\begin{equation}
\label{eq259''}
\inf\limits_{r \in (0,1]}\frac{\mathcal{H}^{d}(Q(x,r) \cap S)}{r^{d}} \geq c_{1}  \quad \hbox{for every} \quad x \in S.
\end{equation}

Then, there exists $\lambda_{0}(S) \in (0,1)$ such that $S$ is $(d,\lambda)$-quasi-porous for every $\lambda \in (0,\lambda_{0}]$.
\end{Lm}

\begin{proof}
The proof is very close in spirit to the proof of Proposition 9.18 in \cite{Tr}. We present the details for the completeness.

Since $\mathcal{H}^{d}(S) \in (0,\infty)$, it is well known (see for example Theorem 1.3.9 in \cite{Lin}) there exists a constant $c_{2} > 0$ depending only on $d$ and there is a set $S' \subset S$ with
$\mathcal{H}^{d}(S \setminus S')=0$ such that
\begin{equation}
\label{eq260''}
\varlimsup\limits_{r \to 0}\frac{\mathcal{H}^{d}(Q(x,r)\cap S)}{r^{d}} \le c_{2} \quad \hbox{for every} \quad x \in S'.
\end{equation}

Now we fix an arbitrary point $x_{0} \in S'$. By \eqref{eq259''} and \eqref{eq260''} there is a
small $r_{0}=r_{0}(x_{0}) > 0$ such that for any $r \in (0,r_{0})$ it holds
\begin{equation}
\label{eq261''}
c_{1} \le \frac{\mathcal{H}^{d}(Q(x_{0},r)\cap S)}{r^{d}} \le 2c_{2}.
\end{equation}
Fix a sufficiently big number $N \in \mathbb{N}$, fix $r \in (0,r_{0}/2)$, subdivide the cube $Q(x_{0},r)$
into $2^{Nn}$ congruent cubes and choose those of them that have a nonempty intersection with $S$. Let
$\{Q_{i}\}_{i=1}^{M}$ be the family of all such chosen cubes. For each $i \in \{1,...,M\}$ take an arbitrary
point $x_{i} \in Q_{i} \cap S$. Clearly, we have $Q_{i} \subset Q(x_{i},2^{-N+1}r):=Q_{i}^{\ast}$ for every $i \in \{1,...,M\}$.
Hence, $S \cap Q(x_{0},r) \subset \cup_{i=1}^{M}Q_{i}^{\ast}$ and the multiplicity of
the covering of $Q(x_{0},r)\cap S$ by the cubes $Q_{i}^{\ast}$, $i=1,...,M$ is bounded above by $5^{n}$. As a result, using
additivity of $\mathcal{H}^{d}$ we get
\begin{equation}
\notag
2c_{2}r^{d} \geq \mathcal{H}^{d}(Q(x_{0},r) \cap S) \geq 5^{-n}\sum\limits_{i=1}^{M}\mathcal{H}^{d}(Q_{i}^{\ast} \cap S)
\geq 5^{-n}c_{1} M 2^{d-Nd} r^{d}.
\end{equation}
Hence, this gives
\begin{equation}
M \le \Bigl(\frac{2 c_{2}}{c_{1}}5^{n}2^{-d}2^{(d-n)N}\Bigr)2^{nN}.
\end{equation}

Since $d < n$  we can take $N=N(d,n,c_{1},c_{2}) \in \mathbb{N}$ big enough to deduce existence of at least one cube $Q_{i} \subset Q(x_{0},r) \setminus S$.
Hence, if we set $\lambda_{0}=1/N$ we get that $x_{0} \in S_{k}(\lambda_{0})$ for every $k > -\log_{2}r+1$. This observation
together with Remark \ref{Rem.porous}
completes the proof.
\end{proof}

\begin{Remark}
Clearly, each Ahlofrs-David $d$-regular set $S \subset \mathbb{R}^{n}$ (with $d \in (0,n)$) satisfies conditions (1) and (2) of Lemma \ref{Lm.porous}.
The converse is false. For example, one can consider a planar rectifiable curve $\Gamma \subset \mathbb{R}^{n}$ with positive length.
\end{Remark}

\subsection{Sobolev spaces}
As usual for each $p \in [1,\infty]$, we let $W_{p}^{1}(\mathbb{R}^{n})$ denote the corresponding Sobolev
space of all equivalence classes of real valued functions $F \in L_{p}(\mathbb{R}^{n})$ whose distributional partial derivatives $D^{\beta}F$ on $\mathbb{R}^{n}$ of order $|\beta| \le 1$ belong to $L_{p}(\mathbb{R}^{n})$.
This space is normed by
\begin{equation}
\notag
\|F\|_{W_{p}^{1}(\mathbb{R}^{n})}:=\sum\limits_{|\beta| \le 1}\|D^{\beta}F\|_{L_{p}(\mathbb{R}^{n})}.
\end{equation}

Recall (see e.g.,\cite{A}, section 6.2) that given a parameter $p \in (1,n]$, for every element $F \in W_{p}^{1}(\mathbb{R}^{n})$ there is \textit{a representative} $\widehat{F}$
in the equivalence class of the element $F$ such that $\widehat{F}$ has
Lebesgue points everywhere except a set $E_{F}$ of $C_{1,p}$-capacity zero. Furthermore, according to imbedding theorem of S. L. Sobolev (see e.g., Theorem 1.2.4 in \cite{A}), given a parameter $p > n$, for every $F \in W_{p}^{1}(\mathbb{R}^{n})$ there is a \textit{continuous representative} $\widehat{F}$ of $F$.
In the sequel we will call $\widehat{F}$ a \textit{good representative
of the element} $F$. Recall also (see Theorem 5.1.13 in \cite{A}) that if $p \in (1,n]$, $d \in (n-p,n]$ then for any given set $S \subset \mathbb{R}^{n}$ the condition $C_{1,p}(S)=0$ implies $\mathcal{H}^{d}(S)=0$. Since in this paper we mainly focus on traces of Sobolev functions to 1-dimensional path-connected sets
in the plane $\mathbb{R}^{2}$ we can present a little bit more rough definition of the trace of a given Sobolev function than the corresponding one
used in \cite{TV}. More precisely, the later appeals to $C_{1,p}$-capacities instead of the Hausdorff measures. However, we believe that the using of $C_{1,p}$-capacities
is not so reasonable in the present framework. Indeed, Proposition \ref{Prop2.5} together with very well known relations between Hausdorff measures
and capacities shows that knowing of a given trace function everywhere except a set of $C_{1,p}$-capacity zero will be overly detailed for us.
Informally speaking, our trace criterion in Theorem 4.2 is expressed in terms of Frostman-type measures $\mu_{k}[\Gamma]$ which "do not
feel" changes of a trace function on a set of $\mathcal{H}^{1}$-measure zero. These remarks justify the following definition.

\begin{Def}
\label{Def2.1}
Let $d \in (0,n]$, $p \in (\max\{1,n-d\},\infty)$ and $F \in W^{1}_{p}(\mathbb{R}^{n})$. Let $S$ be a Borel set with $\mathcal{H}^{d}(S) > 0$.
We define the trace $F|_{S}$ of the element $F$ to the set $S$ as
\begin{equation}
F|_{S}:=\{f: S \to \mathbb{R}: f(x)=\widehat{F}(x) \hbox{ for } \mathcal{H}^{d}-\hbox{a.e. } x \in S\}.
\end{equation}
We define the trace space $W_{p}^{1}(\mathbb{R}^{n})|_{S}$ of the space $W_{p}^{1}(\mathbb{R}^{n})$ as
\begin{equation}
W_{p}^{1}(\mathbb{R}^{n})|_{S}:=\{f: S \to \mathbb{R}: f=F|_{S} \hbox{ for some } F \in W_{p}^{1}(\mathbb{R}^{n})\}
\end{equation}
and equip it with the usual quotient-space norm, i.e.
\begin{equation}
\|f|W_{p}^{1}(\mathbb{R}^{n})|_{S}\|:=\inf\{\|F|W_{p}^{1}(\mathbb{R}^{n})\|: f=F|_{S}\}.
\end{equation}
We denote by $\operatorname{Tr}|_{S}: W_{p}^{1}(\mathbb{R}^{n}) \to W_{p}^{1}(\mathbb{R}^{n})|_{S}$ the corresponding \textit{trace operator.}
\end{Def}

\begin{Remark}
Since the trace $F|_{S}$ of a given Sobolev function is a class  of equivalent (modulo coincidence on a set of $\mathcal{H}^{d}$-measure zero) functions $f:S \to \mathbb{R}$, the trace operator \textit{is well defined  and linear.}
\end{Remark}

\subsection{Rectifiable curves in $\mathbb{R}^{n}$}
By \textit{a curve in $\mathbb{R}^{n}$} we mean the image $\Gamma$ of a continuous map $\gamma: [a,b] \to \mathbb{R}^{n}$
, i.e. $\Gamma=\gamma([a,b])$. The function $\gamma$ is \textit{called a parametrization of the curve} $\Gamma$.
We say that a curve $\Gamma$ is \textit{simple} if the map $\gamma$ is injective.
As usual, we say that a  curve $\Gamma$ is \textit{rectifiable} provided that
\begin{equation}
\label{eq2.1}
l(\Gamma):=\sup \sum\limits_{i=1}^{n}\|\gamma(t_{i})-\gamma(t_{i-1})\| < \infty,
\end{equation}
where the supremum is taken over all tuples $\{t_{i}\}_{i=0}^{n}$ such that
$a=t_{0} < t_{1} < ... < t_{n}=b$.
The associated \textit{length function} $s_{\gamma}$ is defined as
\begin{equation}
\notag
s_{\gamma}(t):=l(\gamma([a,a+t])), \quad t \in [a,b].
\end{equation}

The properties summarised in the next proposition
are well known (see e.g. \cite{Kos}, ch.5).

\begin{Prop}
\label{Prop2.1}
Let $\Gamma \subset \mathbb{R}^{n}$ be a rectifiable curve and let $\gamma: [a,b] \to \mathbb{R}^{n}$
be its parametrization. Then $\Gamma$ admits the $1$-Lipschitz arc length parametrization.
More precisely, there exists the $1$-Lipschitz
map $\gamma_{s}:[0,l(\gamma)] \to \mathbb{R}^{n}$ defined by
\begin{equation}
\label{eq2.2}
\gamma_{s}(\tau):=\gamma(s^{-1}_{\gamma}(\tau)),
\end{equation}
where
\begin{equation}
\label{eq2.3}
s^{-1}_{\gamma}(\tau):=\sup\{s:s_{\gamma}(s)=\tau\}.
\end{equation}
If the curve $\Gamma$ is simple, then
\begin{equation}
\label{eq2.4}
l(\Gamma)=\mathcal{H}^{1}(\Gamma).
\end{equation}
\end{Prop}

It is also useful to recall the infinitesimal behavior of $\mathcal{H}^{1}\lfloor_{\Gamma}$.
The following property is also well known (see e.g. Lemma 3.5 in \cite{Falc})

\begin{Prop}
\label{Prop.curve}
Let $\Gamma \subset \mathbb{R}^{n}$ be a rectifiable curve. Then,
\begin{equation}
\label{eq2.34'}
\lim\limits_{r \to 0}\frac{\mathcal{H}^{1}\lfloor_{\Gamma}(Q(x,r))}{2r}=1 \quad \hbox{for} \quad \mathcal{H}^{1}-\hbox{a.e.} \quad x \in \Gamma.
\end{equation}
\end{Prop}

Clearly, one could assume everywhere in the sequel that any given curve $\Gamma$ has a $1$-Lipschitz parametrization (for example the arc-length parametrization).
However, as we will see in section 5 in practice
it is useful to have a some sort of flexibility in the choice of parameterizations.
This fact justifies the following definition.

\begin{Def}
\label{Def.adm.param}
Given a curve $\Gamma \subset \mathbb{R}^{n}$ we say that a map $\gamma: [a,b] \to \mathbb{R}^{n}$ \textit{is an admissible parametrization of the curve} $\Gamma$
if the following conditions hold:

{\rm (1)} the map $\gamma$ is absolutely continuous;

{\rm (2)}
\begin{equation}
\|\dot{\gamma}(t)\| > 0 \quad \hbox{for} \quad \mathcal{H}^{1}-\hbox{a.e.} \quad t \in [a,b].
\end{equation}

\end{Def}

In the sequel by $\Pi_{i}$, $i=1,...,n$ we denote the \textit{projection maps along} the $i$-th coordinate axes, i.e.
if $x=(x_{1},...,x_{n}) \in \mathbb{R}^{n}$ then $\Pi_{i}(x):=(x_{1},..,x_{i-1},x_{i+1},..,x_{n}) \in \mathbb{R}^{n-1}$. Sometimes we will use the shorthand
$\widehat{x}^{i}:=\Pi_{i}(x)$, $i=1,...,n$ for brevity.
Given a parametrization   $\gamma=(\gamma_{1},...,\gamma_{n}): [a,b] \to \mathbb{R}^{n}$ of the curve $\Gamma$ we let $\gamma_{i}$, $i=1,...,n$ denote
its projections to the $i$-th coordinate axes. Similarly, if for some $t_{0} \in [a,b]$ there exists
the velocity vector
\begin{equation}
\label{eq2.5}
\dot{\gamma}(t_{0}):=\lim\limits_{t \to t_{0}}\frac{\gamma(t)-\gamma(t_{0})}{t-t_{0}} \in \mathbb{R}^{n},
\end{equation}
we let $\dot{\gamma}_{i}(t_{0})$, $i=1,...,n$ denote
the corresponding projections.

\begin{Def}
\label{Def.indicatrix}
Let  $\Gamma \subset \mathbb{R}^{n}$ be a simple rectifiable curve and let $\gamma: [a,b] \to \mathbb{R}^{n}$ be an admissible parametrization of $\Gamma$. 
Given $i \in \{1,...,n\}$, we identify $\mathbb{R}^{n-1}$ with a hyperplane orthogonal to the $i$-th coordinate axe and define
\begin{equation}
\label{eq2.6}
\begin{split}
&\mathcal{L}_{i}[\Gamma](x'):=\Pi_{i}^{-1}(x') \cap \Gamma, \quad x' \in \mathbb{R}^{n-1}.
\end{split}
\end{equation}
Informally speaking $\mathcal{L}_{i}[\Gamma](x')$ is just an intersection of the line going through $x'$ parallel to the $i$-th coordinate axe
with $\Gamma$. We also define the Banach indicatrix functions of $\Gamma$ by
\begin{equation}
\label{eq2.7}
\begin{split}
&N_{i}[\Gamma](x'):=\operatorname{card}\mathcal{L}_{i}[\Gamma](x')=\operatorname{card}\{t \in [a,b]| t=(\Pi_{i} \circ \gamma)^{-1}(x')\}, \quad x' \in \mathbb{R}^{n-1}. \\
\end{split}
\end{equation}
\end{Def}

\begin{Remark}
\label{Rem.measurable}
It follows from Theorem 4.3.2 in \cite{Lin} that for each $i \in \{1,...,n\}$ the function $N_{i}[\Gamma] \in L_{1}^{\operatorname{loc}}(\mathbb{R}^{n-1})$.
Hence, $N_{i}[\Gamma](x') \in \mathbb{N}_{0}$ for each $i \in \{1,...,n\}$ and almost all $x' \in \mathbb{R}^{n-1}$.
Furthermore, analysis of the proof of Lemma 4.1.4 in \cite{Lin} allows to deduce that the functions $N_{i}[\Gamma]$ are Borel.
\end{Remark}

\section{Construction of the special $1$-regular sequence of measures.}

As we mentioned in the introduction, every path-connected set in $\mathbb{R}^{n}$ consisting of more than one point is $1$-thick.
In particular, every  curve $\Gamma \subset \mathbb{R}^{n}$ with $l(\Gamma) > 0$ is $1$-thick. Hence, there is a $1$-regular sequence of measures on $\Gamma$.
The aim of this section is a construction for \textit{any planar simple rectifiable curve} $\Gamma \subset \mathbb{R}^{2}$ \textit{with positive length} of
\textit{a special $1$-regular sequence of measures} $\{\mu_{k}[\Gamma]\}$ concentrated on $\Gamma$. All constructions can be easily extended to
simple rectifiable curves in $\mathbb{R}^{n}$. We consider the case $n=2$ just for the simplicity.

Let us briefly describe the main idea of the construction.
During the section we will assume without loss of generality that $\Gamma \subset [0,1) \times [0,1)$.
By Remark \ref{Rem.measurable} for each $i \in \{1,...,n\}$
intersections of the curve $\Gamma$ with almost every lines $L^{i}$ parallel to the $i$-th coordinate axe consist
of at most finite number of points. Hence, we can easily construct a special $0$-regular sequence of measures concentrated on that finite sets.
After that, taking an arbitrary admissible parametrization $\gamma:[a,b] \to \mathbb{R}^{2}$ of $\Gamma$ and applying Proposition \ref{Prop2.2}
we easily obtain 1-regular sequence of measures concentrated on $\Gamma$.

Let $I_{0}:=[0,1)$ be a half-open unit interval. Given $k \in \mathbb{N}_{0}$, let
\begin{equation}
\notag
\mathcal{D}_{k}:=\Bigl\{I_{k,m}:=\Bigl[\frac{m}{2^{k}},\frac{m+1}{2^{k}}\Bigr), m=0,...,2^{k}-1\Bigr\}
\end{equation}
be the family of all dyadic half-open intervals with side length $2^{-k}$ contained in $I_{0}$.

\begin{Def}
Given a set of distinct points $\{x_{i}\}_{i=0}^{N} \subset I_{0}$, $N \in \mathbb{N}_{0}$, let $\operatorname{P}(\{x_{i}\}_{i=0}^{N})$ be the set of all probability
measures with the support $\{x_{i}\}_{i=0}^{N}$. In other words, $\nu \in \operatorname{P}(\{x_{i}\}_{i=0}^{N})$ if and only if
there exists a \textit{density function} $\alpha: \{x_{i}\}_{i=0}^{N} \to \mathbb{R}_{+}$ with $\sum_{i=0}^{N}\alpha(x_{i})=1$ such that
\begin{equation}
\nu=\sum\limits_{i=0}^{N}\alpha(x_{i})\delta_{x_{i}},
\end{equation}
where $\delta_{x_{i}}$,$i=1,...,N$ are the Dirac measures concentrated at the points $x_{i}$, $i=0,...,N$.
\end{Def}

\begin{Def}
\label{Def.treeforintervals}
Let   $\{x_{i}\}_{i=0}^{N} \subset I_{0}$, $N \in \mathbb{N}$
be  an arbitrary finite set of distinct points. Let $k^{\ast}:=k^{\ast}(\{x_{i}\}_{i=0}^{N})$ be the minimal among all $k \in \mathbb{N}_{0}$ for each of which the map sending
every point $x_{i}$ to the unique dyadic interval $I_{k,m}(x_{i}) \ni x_{i}$ is injective.
For every $k \geq k^{\ast}$ we define the family
\begin{equation}
\mathcal{F}_{k}(\{x_{i}\}_{i=0}^{N}):=\Bigl\{I_{k,m}: I_{k,m} \cap \{x_{i}\}_{i=0}^{N} \neq \emptyset\Bigr\}.
\end{equation}
For each $i\in\{0,...,N\}$ and any $k\geq k^{\ast}(\{x_{i}\}_{i=0}^{N})$ we denote by $I_{k,m}(x_{i})$ the unique dyadic interval
in $\mathcal{F}_{k}(\{x_{i}\}_{i=0}^{N})$ containing the point $x_{i}$.
\end{Def}

Given an arbitrary tree $\mathcal{T}$ with a root $r$ we introduce the intrinsic metric $\rho$ on  $\mathcal{T}$ making it the so-called metric tree
$(\mathcal{T},\rho)$.
More precisely, given two vertices $\xi,\xi' \in \mathcal{V}(\mathcal{T})$ joined by some edge $e$ we put $\rho(\xi,\xi')=1$. For generic two
vertices $\xi,\xi' \in \mathcal{V}(\mathcal{T})$  we define
$\rho(\xi,\xi')=\inf\sum_{i}\rho(\xi_{i},\xi_{i-1})$, where the infimum is taken over all paths
$\xi=:\xi_{0} \leftrightarrow .... \leftrightarrow \xi_{l}:=\xi'$ joining $\xi$ with $\xi'$.
Given $i \in \mathbb{N}_{0}$, let $\mathcal{V}^{i}(\mathcal{T}):=\{\xi \in \mathcal{T}:\rho(r,\xi)=i\}$.
If $\mathcal{V}^{i}(\mathcal{T}) \neq \emptyset$ then for any $\xi \in \mathcal{V}^{i}(\mathcal{T})$ we denote by $n(\xi)$
the number of edges joining $\xi$ with the corresponding vertices in $\mathcal{V}^{i+1}(\mathcal{T})$.

\begin{Def}
Let $k \in \mathbb{N}_{0}$. Given a nonempty family of dyadic intervals $\mathcal{F}_{k} \subset \mathcal{D}_{k}$, we define
the tree $\mathcal{T}=\mathcal{T}(\mathcal{F}_{k})$ as the metric tree whose vertices $\mathcal{V}^{i}(\mathcal{T})$, $i=0,...,k$
are naturally corresponds to all those dyadic intervals in $\mathcal{D}_{i}$ each of which contains at least one interval from $\mathcal{F}_{k}$.
\end{Def}

\begin{Def}
\label{Def.Spec.Frostman}
Let $\mathcal{F}_{k} \subset \mathcal{D}_{k}$, $k \in \mathbb{N}_{0}$ be an arbitrary nonempty
family of dyadic intervals. For each $j \in \mathbb{N}_{0}$ and any $I \in \mathcal{F}_{k}$  we define \textit{the Frostman-type weight} by letting
\begin{equation}
\label{eq2.9}
\alpha^{j}_{F}[\mathcal{F}_{k}](I):=
\begin{cases}
&\prod\limits_{l=j}^{k-1}\frac{1}{n(\xi_{l})} \quad \hbox{if} \quad j \in \{0,...,k-1\} ,\\
&1 \quad \hbox{if} \quad j > k-1,
\end{cases}
\end{equation}
where $\xi(I) \in \mathcal{V}^{k}(\mathcal{T})$ is the unique vertex corresponding to the interval $I$
and $r=:\xi_{0} \leftrightarrow .... \leftrightarrow \xi_{k}:=\xi(I)$ is the unique path joining the root $r$
with the vertex $\xi(I)$.
Given a finite set $\{x_{i}\}_{i=0}^{N}$, $N \in \mathbb{N}$ of distinct points, we define \textit{the Frostman-type weight} by letting
\begin{equation}
\label{eq2.9''}
\alpha^{j}_{F}[\{x_{i}\}_{i=0}^{N}](x_{i}):=
\alpha^{j}_{F}[\mathcal{F}_{k}(\{x_{i}\}_{i=0}^{N})](I_{k,m}(x_{i})), \quad k \geq k^{\ast}(\{x_{i}\}_{i=0}^{N}).
\end{equation}
\end{Def}

\begin{Def}
Given a finite set $\{x_{i}\}_{i=0}^{N} \subset I_{0}$, $N \in \mathbb{N}$ of distinct points, we define for each $j \in \mathbb{N}_{0}$ \textit{the Frostman-type  measure} $\nu^{j}_{F}[\{x_{i}\}_{i=0}^{N}]$ by letting
\begin{equation}
\nu^{j}_{F}[\{x_{i}\}_{i=0}^{N}]:=\sum\limits_{i=0}^{N}\alpha^{j}_{F}[\{x_{i}\}_{i=0}^{N}](x_{i})\delta_{x_{i}}.
\end{equation}
\end{Def}

\begin{Remark}
It follows immediately from the definition of $k^{\ast}$ that Definition \ref{Def.Spec.Frostman} is correct.
Indeed, for each $k \geq k^{\ast}$ we have $\alpha^{j}_{F}[\mathcal{F}_{k}(\{x_{i}\}_{i=0}^{N})](I)=
\alpha^{j}_{F}[\mathcal{F}_{k^{\ast}}(\{x_{i}\}_{i=0}^{N})](I)$ for all $I \in \mathcal{F}_{k}(\{x_{i}\}_{i=0}^{N})$.
\end{Remark}

The following assertion exhibits basic properties of Frostman-type weights.

\begin{Lm}
\label{Lm2.1}
Let $\{x_{i}\}_{i=0}^{N} \subset I_{0}$, $N \in \mathbb{N}_{0}$
be a set of distinct points. Then, for each $j \in \mathbb{N}_{0}$ the following holds:

{\rm (1)} for and any dyadic interval
$I_{j,m}$, $m \in \{0,...,2^{j}-1\}$
\begin{equation}
\label{eq211}
\nu^{j}_{F}[\{x_{i}\}_{i=0}^{N}](I_{j,m}) =
\begin{cases}
&1, \quad \{x_{i}\}_{i=0}^{N} \cap I_{j,m} \neq \emptyset;\\
&0, \quad \{x_{i}\}_{i=0}^{N} \cap I_{j,m} = \emptyset.
\end{cases}
\end{equation}

{\rm (2)} for every $i \in \{0,...,N\}$ either
\begin{equation}
\begin{split}
\label{eq212}
&\alpha^{j}_{F}[\{x_{i}\}_{i=0}^{N}](x_{i}) = \alpha^{j+1}_{F}[\{x_{i}\}_{i=0}^{N}](x_{i}) \quad \hbox{or}\\
&\alpha^{j}_{F}[\{x_{i}\}_{i=0}^{N}](x_{i})= 2^{-1}\alpha^{j+1}_{F}[\{x_{i}\}_{i=0}^{N}](x_{i}).
\end{split}
\end{equation}
\end{Lm}

\begin{proof}
For $j \geq k^{\ast}$ the statement is obvious. Let $\mathcal{T}=\mathcal{T}(\mathcal{F}_{k^{\ast}}(\{x_{i}\}_{i=0}^{N}))$ be the corresponding metric tree. Note that in the case $j < k^{\ast}$ equalities \eqref{eq212} follow directly from \eqref{eq2.9} because for any $\xi \in \mathcal{V}(\mathcal{T})$ we have either $n(\xi)=1$ or $n(\xi)=2$.

To prove \eqref{eq211} we argue by induction. For $j \geq k^{\ast}$ this is obvious. Suppose that $k^{\ast} > 1$ and that \eqref{eq211} is proved for some $j_{0} \in \{1,...,k^{\ast}\}$.
Then from \eqref{eq212} it is easy to conclude the validity of \eqref{eq211} for $j_{0}-1$ using the same arguments as above.
\end{proof}

\begin{Remark}
\label{Remm32}
Note that given a set $\{x_{i}\}_{i=0}^{N} \subset I_{0}$, $N \in \mathbb{N}_{0}$
of distinct points and $I_{j,m}$ with $I_{j,m} \cap \{x_{i}\}_{i=0}^{N} \neq \emptyset$,
the restriction $\nu^{j}_{F}\lfloor_{I_{j,m}} \in \operatorname{P}(I_{j,m} \cap \{x_{i}\}_{i=0}^{N})$. Of course in general
$\nu^{j}_{F}$ is not a probability measure on $I_{0}$.
\end{Remark}

\begin{Def}
\label{Def2.3}
Let $\Gamma \subset [0,1) \times [0,1)$ be a simple rectifiable curve with positive length. Let $\gamma: [a,b] \to \mathbb{R}^{2}$ be an admissible parametrization of $\Gamma$.
 Given $j \in \mathbb{N}_{0}$, we define for each $t \in [a,b]$
\begin{equation}
\label{eq2.11}
\mathcal{W}^{j}_{1}[\gamma](t):=
\begin{cases}
&\alpha^{j}_{F}[\mathcal{L}_{1}[\Gamma](\gamma_{2}(t))](\gamma_{1}(t)) \quad \hbox{if} \quad \quad N_{1}[\Gamma](\gamma_{2}(t)) < \infty;\\
&0 \quad \hbox{if} \quad \quad N_{1}[\Gamma](\gamma_{2}(t)) = \infty;
\end{cases}
\end{equation}
and similarly,
\begin{equation}
\label{eq2.12}
\mathcal{W}^{j}_{2}[\gamma](t):=
\begin{cases}
&\alpha^{j}_{F}[\mathcal{L}_{2}[\Gamma](\gamma_{1}(t))](\gamma_{2}(t)) \quad \hbox{if} \quad \quad N_{2}[\Gamma](\gamma_{1}(t)) < \infty;\\
&0 \quad \hbox{if} \quad \quad N_{2}[\Gamma](\gamma_{1}(t)) = \infty.
\end{cases}
\end{equation}
Given  $j \in \mathbb{N}_{0}$ we also define \textit{the special density} as
\begin{equation}
\label{spec.density}
\operatorname{D}^{j}[\gamma](t):=\max\{\mathcal{W}^{j}_{2}[\gamma](t)|\dot{\gamma}_{1}(t)|,\mathcal{W}^{j}_{1}[\gamma](t)|\dot{\gamma}_{2}(t)|\},
\quad t \in [a,b].
\end{equation}
\end{Def}

\begin{Remark}
\label{Remm3.3}
Let $\gamma: [a,b] \to \mathbb{R}^{2}$ and $\gamma': [a',b'] \to \mathbb{R}^{2}$ be two admissible parameterizations of $\Gamma$.
Suppose that $\gamma(t)=\gamma'(t'(t))$ for some strictly increasing absolutely continuous function $t':[a,b] \to [a',b']$.
Then it is easy to show that $\mathcal{W}^{j}_{i}[\gamma](t)=\mathcal{W}^{j}_{i}[\gamma'](t'(t))$ for each $i=1,2$, every $j \in \mathbb{N}_{0}$
and all $t \in [a,b]$.
\end{Remark}

\begin{Prop}
Let $\Gamma \subset [0,1) \times [0,1)$ be a simple rectifiable curve with positive length. Let $\gamma: [a,b] \to \mathbb{R}^{2}$ be an admissible parametrization
of $\Gamma$.
Then for every $j \in \mathbb{N}_{0}$  the functions $\mathcal{W}^{j}_{i}[\gamma]$, $i=1,2$ and $\operatorname{D}^{j}[\gamma]$ are measurable.
\end{Prop}

\begin{proof}
First of all we prove that the functions $\mathcal{W}^{j}_{1}[\gamma]$, $j \in \mathbb{N}_{0}$ are measurable.
The proof of measurability of the functions $\mathcal{W}^{j}_{2}[\gamma]$, $j \in \mathbb{N}_{0}$ requires exactly the same arguments.
By Lemma \ref{Lm2.1} it is sufficient to verify that $\mathcal{W}^{0}_{1}[\gamma]$ is measurable.

Fix an arbitrary number $k \in \mathbb{N}$. Let $\mathcal{F}_{k} \subset \mathcal{D}_{k}$ be an arbitrary nonempty family of dyadic intervals.
We define
\begin{equation}
\notag
\begin{split}
&F_{k}:=\bigcup\limits_{I_{k,m} \in \mathcal{F}_{k}}I_{k,m}, \quad \widetilde{E}(F_{k}):=\gamma_{2}(\gamma_{1}^{-1}(F_{k}))\setminus \gamma_{2}([a,b]\setminus \gamma_{1}^{-1}(F_{k})).
\end{split}
\end{equation}
In other words, $\widetilde{E}(F_{k})$ is the set of all $x_{2} \in \Pi_{1}(\Gamma)$ such that $\mathcal{L}_{1}[\Gamma](x_{2}) \subset F_{k}$.
The set $\gamma_{1}^{-1}(F_{k})$ is Borel. Since $\gamma_{2}$ is continuous, this implies that the set
$\widetilde{E}(F_{k})$ is a difference of two Souslin sets and
hence, it is universally measurable.
Now we define
\begin{equation}
\notag
E(F_{k}):=\widetilde{E}(F_{k}) \setminus \bigcup\limits_{\substack{F'_{k} \subset F_{k} \\ F'_{k} \neq F_{k}}}\widetilde{E}(F'_{k}),
\quad G(F_{k}):=\gamma_{2}^{-1}(E(F_{k})).
\end{equation}
Informally speaking, $E(F_{k})$ is a set of all points on the second coordinate axe such that
for all lines going trough that points and parallel to the first coordinate axe the corresponding intersections
with $\Gamma$ consist of the sets of points whose projections to the first coordinate axe meet every interval
from the family $\mathcal{F}_{k}$ and do not meet the other dyadic intervals.

Clearly, $E(F_{k})$ and $G(F_{k})$ are universally measurable. We set
\begin{equation}
\mathcal{W}_{k,1}[\gamma](t):=\alpha_{F}[\mathcal{F}_{k}](I_{k,m}), \quad \hbox{if} \quad t \in \gamma^{-1}_{1}(I_{k,m}) \cap G(F_{k}).
\end{equation}
Since for different families $\mathcal{F}_{k}$ and $\mathcal{F}'_{k}$ the sets $G(F_{k})$
and $G(F'_{k})$ are disjoint and since
\begin{equation}
\notag
[a,b]=\bigcup\limits_{\mathcal{F}_{k}} G(F_{k})
\end{equation}
the function $\mathcal{W}_{k,1}[\gamma]$ is well defined everywhere on $[a,b]$ and measurable.

From Remark \ref{Rem.measurable} it follows that  the set $\mathcal{G}$ where $N_{1}[\Gamma]=+\infty$ has measure $\mathcal{H}^{1}(\mathcal{G})=0$.
Since the function $N_{1}[\Gamma]$ is Borel, the set $\mathcal{G}$ is Borel. Hence, the set $\gamma^{-1}_{2}(\mathcal{G})$ is Borel.
As a result, by the very definition of the sets $E(F_{k})$, $k \in \mathbb{N}_{0}$ and $\mathcal{G}$ it follows that
\begin{equation}
\notag
\lim\limits_{k \to \infty}\mathcal{W}_{k,1}[\gamma](t)=\mathcal{W}^{0}_{1}[\gamma](t) \quad \hbox{for every} \quad t \in [a,b] \setminus
\gamma^{-1}_{2}(\mathcal{G}).
\end{equation}
This implies that  $\mathcal{W}^{0}_{1}[\gamma]$ is measurable.

Note that the image of the function $\mathcal{W}^{0}_{1}[\gamma]$ is an at most countably set. Using this fact together with the measurability of
$\mathcal{W}^{0}_{1}[\gamma]$ and \eqref{eq2.9}, \eqref{eq2.9''}, \eqref{eq2.11} it is easy to get measurability of $\mathcal{W}^{j}_{1}[\gamma]$
for all $j \in \mathbb{N}$.

Finally, to prove measurability of $\operatorname{D}^{j}[\gamma]$ it sufficient to use measurability of $\mathcal{W}^{j}_{i}$ just established
and take into account item (1) of Definition \ref{Def.adm.param}.

The proof is complete.

\end{proof}

Now we are ready to define the notion which will be the keystone in this paper.

\begin{Def}
\label{Def2.4}
Let $\Gamma \subset  [0,1) \times [0,1)$ be a simple rectifiable curve with positive length. Let $\gamma: [a,b] \to \mathbb{R}^{2}$ be an admissible parametrization of $\Gamma$.
We say that the sequence of measures $\{\mu_{k}[\Gamma]\}:=\{\mu_{k}[\Gamma]\}_{k \in \mathbb{N}_{0}}$ \textit{is a special sequence of measures on} $\Gamma$ if
\begin{equation}
\label{eq2.13}
\mu_{k}[\Gamma]:=\gamma_{\sharp}(\operatorname{D}^{k}[\gamma]\mathcal{H}^{1}\lfloor_{[a,b]}) \quad \hbox{for every} \quad k \in \mathbb{N}_{0}.
\end{equation}
\end{Def}

We are going to show that the measures $\mu_{k}[\Gamma]$, $k \in \mathbb{N}_{0}$ are well defined.
For that purpose we recall the following property. Probably it looks like a folklore but we present the proof for the completeness.

\begin{Prop}
\label{Propp.3.2}
Let $[a,b] \subset \mathbb{R}$, $[A,B] \subset \mathbb{R}$ and $g: [a,b] \to [A,B]$ be absolutely continuous. Assume that
\begin{equation}
\label{eqq.313}
\frac{dg}{dt}(t) > 0 \quad \hbox{for} \quad \mathcal{H}^{1}-a.e. \quad  t \in [a,b].
\end{equation}
Then the inverse function $g^{-1}: [A,B] \to [a,b]$ is absolutely continuous.
\end{Prop}

\begin{proof}
From \eqref{eqq.313} it follows that $g: [a,b] \to [A,B]$ is strictly increasing.
Therefore, there exists the inverse function $g^{-1}: [A,B] \to [a,b]$ which is continuous and strictly increasing. Hence, to prove that
$g^{-1}: [A,B] \to [a,b]$ is absolutely continuous it is sufficient to check the Lusin property.
Assume the contrary. Hence, there exists a set $E \subset [A,B]$ with $\mathcal{H}^{1}(E)=0$ such that $\mathcal{H}^{1}(g^{-1}(E)) > 0$.
Combing this with \eqref{eqq.313} we get
\begin{equation}
\label{eqq.314}
r_{0}:=\int\limits_{g^{-1}(E)}\frac{dg}{dt}(t) dt > 0.
\end{equation}
Let $\{(A_{i},B_{i})\}_{i=1}^{\infty}$ be an arbitrary sequence of nonempty intervals such that $E \subset \cup_{i=1}^{\infty}(A_{i},B_{i})$.
Since $g^{-1}$ is strictly increasing we have a sequence of nonempty intervals $\{(a_{i},b_{i})\}_{i=1}^{\infty}$ such that $(a_{i},b_{i}):=g^{-1}((A_{i},B_{i}))$, $i \in \mathbb{N}$. Clearly $g^{-1}(E) \subset \cup_{i=1}^{\infty}(a_{i},b_{i})$.
Since $g$ is absolutely continuous we can apply Newton-Leibniz formula  and take into account \eqref{eqq.314}. As a result, we have
\begin{equation}
\label{eqq.315}
\sum\limits_{i=1}^{\infty}|A_{i}-B_{i}|= \sum\limits_{i=1}^{\infty} \int\limits_{a_{i}}^{b_{i}}\frac{dg}{dt}(t) dt \geq r_{0}.
\end{equation}
The sequence $\{(A_{i},B_{i})\}_{i=1}^{\infty}$ was chosen arbitrarily. Hence, by \eqref{eqq.315} and the definition of the measure $\mathcal{H}^{1}$ we conclude
that $\mathcal{H}^{1}(E) \geq r_{0} > 0$.
This contradiction completes the proof.
\end{proof}

\begin{Prop}
\label{Prop.indep.parametr.}
Let $\Gamma\subset [0,1) \times [0,1)$ be a simple rectifiable curve with positive length. Let $\gamma: [a,b] \to \mathbb{R}^{2}$ and $\gamma': [a',b']
 \to \mathbb{R}^{2}$ be  admissible parameterizations of $\Gamma$.
Then,
\begin{equation}
\label{eq.indep.parametr.}
\gamma_{\sharp}(\operatorname{D}^{k}[\gamma]\mathcal{H}^{1}\lfloor_{[a,b]})=\gamma'_{\sharp}(\operatorname{D}^{k}[\gamma']\mathcal{H}^{1}\lfloor_{[a',b']}).
\end{equation}
\end{Prop}

\begin{proof}
Since the parameterizations $\gamma$ and $\gamma'$ are admissible we obtain that the associated length functions $s_{\gamma}$ and
$s_{\gamma'}$ are absolutely continuous and strictly increasing.
Furthermore,
\begin{equation}
\label{eq3.16''''}
\begin{split}
&\dot{s}_{\gamma}(t) = \|\dot{\gamma}(t)\|> 0 \quad \hbox{for}  \quad \mathcal{H}^{1}-\hbox{a.e.} \quad t \in [a,b];\\
&\dot{s}_{\gamma'}(t') = \|\dot{\gamma'}(t')\|> 0 \quad \hbox{for}  \quad \mathcal{H}^{1}-\hbox{a.e.} \quad t' \in [a',b'].
\end{split}
\end{equation}
Hence, there are continuous strictly increasing functions $t=t(s)=s_{\gamma}^{-1}(s)$ and $t'=t'(s)=s_{\gamma'}^{-1}(s)$. We set
$t'(t):=  s_{\gamma'}^{-1}(s_{\gamma}(t))$ for all $t \in [a,b]$. Note that
\begin{equation}
\notag
\frac{dt'}{dt}(t) > 0, \quad \hbox{for} \quad \mathcal{H}^{1}-a.e. \quad t \in [a,b].
\end{equation}
By Proposition \ref{Propp.3.2} not only the functions $s_{\gamma}$, $s_{\gamma'}$ but also the functions $s_{\gamma}^{-1}$, $s_{\gamma'}^{-1}$ have the
Lusin property. Hence, the function $t'(\cdot)$ is absolutely continuous and strictly increasing.
As a result, we apply Proposition \ref{Prop2.2} with $\Phi(t)=t'(t)$ taking into account \eqref{spec.density} and Remark \ref{Remm3.3}. This allows to establish
for any measurable set $E$ the desirable equality
\begin{equation}
\notag
\int\limits_{E}\operatorname{D}^{k}[\gamma](t)\,dt=\int\limits_{E}\operatorname{D}^{k}[\gamma'](t'(t))\frac{dt'}{dt}(t)\,dt = \int\limits_{t'(E)}\operatorname{D}^{k}[\gamma'](t')\,dt'.
\end{equation}

The proof is complete.
\end{proof}

Recall Proposition \ref{Prop2.1}.
Now we can formulate \textit{the main result of this section.}

\begin{Th}
\label{Th2.1}
Let $\Gamma \subset [0,1) \times [0,1)$ be a simple rectifiable curve with positive length.
Then the special sequence of measures $\{\mu_{k}[\Gamma]\}:=\{\mu_{k}[\Gamma]\}_{k \in \mathbb{N}_{0}}$ is 1-regular on $\Gamma$.
\end{Th}

\begin{proof}
Using Proposition \ref{Prop.indep.parametr.} we may assume without loss of generality that
the curve $\Gamma$ is parameterized by the arc length $\gamma_{s}$. During the proof we use the shorthand
$\gamma=\gamma_{s}$. By Proposition \ref{Prop2.1} the map $\gamma:[0,l(\Gamma)] \to \Gamma$ is $1$-Lipschitz.
Clearly the maps $\gamma_{i}$, $i=1,2$ are $1$-Lipschitz as compositions of
$\gamma$ with the corresponding projections $\Pi_{i}$.
We set
\begin{equation}
\notag
l_{i}=\operatorname{diam}\gamma_{i}([a,b]) \quad \hbox{for} \quad i=1,2.
\end{equation}
Since the curve $\Gamma$ has a positive length it holds
\begin{equation}
0 < l(\Gamma) \le l_{1}+l_{2} \le 2l(\Gamma) < +\infty.
\end{equation}
Now we should verify that conditions (1)--(4) in the definition of $1$-regular on $\Gamma$ sequence of measures (see the introduction)
are holds true for the sequence $\{\mu_{k}[\Gamma]\}$.

\textit{Step 1.} Since $\gamma$ is the arc length parametrization of $\Gamma$ it follows directly from \eqref{spec.density} that
\begin{equation}
\notag
\operatorname{D}^{k}[\gamma](t) > 0 \quad \hbox{for} \quad \mathcal{H}^{1}-\hbox{a.e.} \quad t \in [a,b].
\end{equation}
Hence, the construction \eqref{eq2.13} gives
\begin{equation}
\label{eq312}
\operatorname{supp}\mu_{k}[\Gamma]=\Gamma \quad \hbox{for every} \quad k \in \mathbb{N}_{0}.
\end{equation}

\textit{Step 2.} Now it is convenient to verify condition (4). It follows immediately from \eqref{eq212}, \eqref{eq2.11} and \eqref{eq2.12} that for every
$k \in \mathbb{N}_{0}$ it holds
\begin{equation}
\label{eq313}
2^{-1}\operatorname{D}^{k+1}[\gamma](t) \le \operatorname{D}^{k}[\gamma](t) \le \operatorname{D}^{k+1}[\gamma](t) \quad \hbox{for} \quad \mathcal{H}^{1}-\hbox{a.e.} \quad t \in [a,b].
\end{equation}
This and \eqref{spec.density}, \eqref{eq2.13} clearly imply existence of a sequence of weights $\{w_{k}[\Gamma]\}$ such that $\mu_{k}[\Gamma]=w_{k}[\Gamma]\mu_{0}[\Gamma]$ for every
$k \in \mathbb{N}$ and
\begin{equation}
\label{eq314}
2^{-1}w_{k+1}[\Gamma](x) \le w_{k}[\Gamma](x) \le w_{k+1}[\Gamma](x) \quad \hbox{for} \quad \mu_{0}[\Gamma]-\hbox{a.e.} \quad x \in \Gamma.
\end{equation}

\textit{Step 3.}
Fix an arbitrary $j \in \mathbb{N}_{0}$. If $i=1$ we set $i'=2$, if $i=2$ we set $i'=1$.
It is clear that for any (half-open) dyadic cube $Q_{j,m}$ with $j \geq k$, $m \in \mathbb{Z}^{2}$  we have by \eqref{spec.density} and \eqref{eq2.13}
\begin{equation}
\begin{split}
\label{eq2.20}
&\int\limits_{\gamma^{-1}(Q_{j,m})}\mathcal{W}^{k}_{i'}[\gamma](s)|\dot{\gamma}_{i}(s)|\,ds
\le \mu_{k}[\Gamma](Q_{j,m})\\
&\le \int\limits_{\gamma^{-1}(Q_{j,m})}
\Bigl(\mathcal{W}^{k}_{1}[\gamma](s)|\dot{\gamma}_{2}(s)|+\mathcal{W}^{k}_{2}[\gamma](s)|\dot{\gamma}_{1}(s)|\Bigr)\,ds, \quad i=1,2.
\end{split}
\end{equation}
Since the map $\gamma$ is injective we have for $\mathcal{H}^{1}$-a.e. $x_{i} \in \Pi_{i'}(Q_{j,m})$ equality
\begin{equation}
\label{eqq323}
\int\limits_{\gamma^{-1}(Q_{j,m}) \cap \gamma^{-1}_{i}(x_{i})} \mathcal{W}^{k}_{i'}[\gamma](s)\,d\mathcal{H}^{0}(s) = \int\limits_{
\Pi_{i}(Q_{j,m})}
d\nu_{F}^{k}[\mathcal{L}_{i'}[\Gamma](x_{i})](x_{i'}).
\end{equation}
Clearly, $\gamma$, $i=1,2$ are Lipschitz. Hence, we apply Proposition \ref{Prop2.2} and use \eqref{eqq323} taking into account Definition \ref{Def2.3} and Remark \ref{Remm32}. 
We get for  each $k \in \mathbb{N}_{0}$ and any $j \geq k$, $m \in \mathbb{Z}^{2}$,
\begin{equation}
\label{eq2.21}
\begin{split}
&\int\limits_{\gamma^{-1}(Q_{j,m})}
\mathcal{W}^{k}_{i'}[\gamma](s)|\dot{\gamma}_{i}(s)|\,ds = \int\limits_{\Pi_{i'}(Q_{j,m} \cap \Gamma)}\Bigl(\int\limits_{
\Pi_{i}(Q_{j,m})}
d\nu_{F}^{k}[\mathcal{L}_{i'}[\Gamma](x_{i})](x_{i'})\Bigr)\,d\mathcal{H}^{1}(x_{i})\\
&\le \int\limits_{\Pi_{i'}(Q_{j,m} \cap \Gamma)}\,d\mathcal{H}^{1}(y)=\mathcal{H}^{1}(\Pi_{i'}(Q_{j,m} \cap \Gamma)), \quad i=1,2.
\end{split}
\end{equation}
On the other hand, similar arguments allow to deduce for every $k \in \mathbb{N}_{0}$, $m \in \mathbb{Z}^{2}$ the following equality
\begin{equation}
\label{eq2.21'}
\int\limits_{\gamma^{-1}(Q_{k,m})}
\mathcal{W}^{k}_{i'}[\gamma](s)|\dot{\gamma}_{i}(s)|\,ds  = \mathcal{H}^{1}(\Pi_{i'}(Q_{k,m} \cap \Gamma)), \quad i=1,2.
\end{equation}

\textit{Step 4.} We verify condition (3) for the sequence $\{\mu_{k}[\Gamma]\}_{k \in \mathbb{N}_{0}}$.
Fix $k \in \mathbb{N}_{0}$ and $x \in \Gamma$. We set $Q_{k}(x):=Q(x,2^{-k})$ for brevity.
We use \eqref{eq312}, \eqref{eq314}, \eqref{eq2.21'} and subadditivity of the measure $\mathcal{H}^{1}$. This gives
\begin{equation}
\label{eq2.23}
\begin{split}
&\mu_{k}[\Gamma](Q_{k}(x) \cap \Gamma)=\mu_{k}[\Gamma](Q_{k}(x)) \geq \frac{1}{4}\mu_{k+2}[\Gamma](Q_{k}(x)) \geq \frac{1}{4}\sum\limits_{\substack{m \in \mathbb{Z}^{2} \\ Q_{k+2,m} \cap \frac{1}{2}Q_{k}(x) \neq \emptyset}}\mu_{k+2}[\Gamma](Q_{k+2,m})\\
&\geq \sum\limits_{\substack{m \in \mathbb{Z}^{2} \\ Q_{k+2,m} \cap \frac{1}{2}Q_{k}(x) \neq \emptyset}}\frac{1}{4}\mathcal{H}^{1}(\Pi_{i}(Q_{k+2,m} \cap \Gamma))
\geq \frac{1}{4}\mathcal{H}^{1}(\Pi_{i}(\frac{1}{2}Q_{k}(x) \cap \Gamma)), \quad i=1,2.
\end{split}
\end{equation}
Since $Q_{k}(x) \cap \Gamma$ is path-connected we clearly have for $k \geq -\log_{2}(\max\{l_{1},l_{2}\})$
\begin{equation}
\label{eq2.24}
\max\{\mathcal{H}^{1}(\Pi_{1}(\frac{1}{2}Q_{k}(x) \cap \Gamma)),
\mathcal{H}^{1}(\Pi_{2}(\frac{1}{2}Q_{k}(x) \cap \Gamma))\} \geq 2^{-k-1}.
\end{equation}
Finally, combining \eqref{eq2.23} and \eqref{eq2.24} we deduce
\begin{equation}
\label{eq2.25}
\mu_{k}[\Gamma](Q_{k}(x)) \geq \frac{1}{8}\min\{1,2^{k}\max\{l_{1},l_{2}\}\}2^{-k} \geq \frac{1}{8}\min\{1,\max\{l_{1},l_{2}\}\}2^{-k}, \quad k \in \mathbb{N}_{0}.
\end{equation}

\textit{Step 5.}
To verify condition (2) we fix $r \in (0,2^{-k}]$ and set $k(r):=[\log_{2}r^{-1}]$. Note that there are at most $25$ dyadic cubes $Q_{k(r),m}$ whose intersections with $Q=Q(x,r)$
are nonempty. Hence, we apply the second inequality in \eqref{eq2.20} and then \eqref{eq2.21} with $j=k(r)$ (it is possible because $k(r) \geq k$). We get
\begin{equation}
\begin{split}
&\mu_{k}[\Gamma](Q) \le \sum\limits_{\substack{m \in \mathbb{Z}^{2} \\ Q_{k(r),m} \cap Q \neq \emptyset}}\mu_{k}[\Gamma](Q_{k(r),m})\\
&\le \sum\limits_{\substack{m \in \mathbb{Z}^{n} \\ Q_{k(r),m} \cap Q \neq \emptyset}}
\sum\limits_{i=1}^{2}\mathcal{H}^{1}(\Pi_{i}(Q_{k(r),m} \cap \Gamma)) \le \frac{50}{2^{k(r)}} \le 50r.
\end{split}
\end{equation}

The proof is complete.

\end{proof}

\section{Main results}

During the whole section we use the shorthand  $Q_{k}(x):=Q(x,2^{-k})$, $k \in \mathbb{N}_{0}$.

The proof of the following lemma is based on standard arguments. Nevertheless, as far as we know, the assertion is new.
We present the full proof for the completeness.

\begin{Lm}
\label{Prop2.6}
Let $\Gamma \subset [0,1) \times [0,1)$ be a simple rectifiable curve of positive length. Let $\{\mu_{k}[\Gamma]\}$ be the special $1$-regular sequence of measures on $\Gamma$. Let $\mu_{k}[\Gamma]=w_{k}[\Gamma]\mu_{0}[\Gamma]$, $k\in \mathbb{N}$.
Then, 
\begin{equation}
\varlimsup\limits_{k \to \infty}\frac{1}{w_{k}[\Gamma](x)}\fint\limits_{Q(x,2^{-k})}w_{k}[\Gamma](y)\,d\mu_{0}[\Gamma](y) < +\infty
\quad \hbox{for} \quad \mathcal{H}^{1}-a.e. \quad x \in \Gamma.
\end{equation}
\end{Lm}

\begin{proof}  Due to Proposition \ref{Prop.indep.parametr.} we may assume that $\Gamma$
is parameterized by the arc length $\gamma_{s}:[0,l(\Gamma)] \to \mathbb{R}^{2}$. We split the proof into several steps.

\textit{Step 1.} Since $\mathcal{H}^{1}\lfloor_{\Gamma}$ is a Radon measure, by Theorem \ref{Th2.1} and Proposition \ref{Prop2.5}
there exists a function $g[\Gamma] \in L_{1}(\mathcal{H}^{1}\lfloor_{\Gamma})$
such that for $\mathcal{H}^{1}$-a.e. point $x \in \Gamma$ it holds (note that $w_{0}[\Gamma] \equiv 1$)
\begin{equation}
\label{eq.41}
g[\Gamma](x)=\frac{d\mu_{0}[\Gamma]}{d\mathcal{H}^{1}\lfloor_{\Gamma}}(x) \le C^{1}_{\{\mu_{k}[\Gamma]\}}.
\end{equation}

\textit{Step 2.} At this step we are going to show that
\begin{equation}
\label{eq264''}
g[\Gamma](x) > 0 \quad \hbox{for}  \quad \mathcal{H}^{1}-\hbox{a.e.} \quad x \in \Gamma.
\end{equation}

By Proposition \ref{Prop.curve} it follows from \eqref{eq.41} that in order to prove \eqref{eq264''} it is sufficient to establish
\begin{equation}
\varliminf\limits_{r \to 0}\frac{\mu_{0}[\Gamma](Q(x,r))}{r} > 0 \quad \hbox{for}  \quad \mathcal{H}^{1}-\hbox{a.e.} \quad x \in \Gamma.
\end{equation}

Since the set $\Gamma$ is connected and the map $\gamma_{s}$ is the arc-length parametrisation of $\Gamma$ 
we get for each point $x \in \Gamma$ and any $r < \frac{\operatorname{diam}\Gamma}{3}$
that the preimage $\gamma^{-1}_{s}(Q(x,r))$
contains a closed interval $[t^{1}_{x}(r),t^{2}_{x}(r)] \ni \gamma^{-1}_{s}(x)$ such that:

{\rm (1)} $\gamma_{s}(t^{1}_{x}(r))=x$, $\gamma_{s}(t^{2}_{x}(r)) \in \partial Q(x,r)$;

{\rm (2)} it holds
\begin{equation}
\label{eq4.6}
|t^{1}_{x}(r)-t^{2}_{x}(r)| \geq r.
\end{equation}
But then, since $\gamma_{s}$ is the arc-length parametrization, taking into account Proposition \ref{Prop.curve} we obtain
\begin{equation}
\label{eq4.7}
|t^{1}_{x}(r)-t^{2}_{x}(r)| \le l(\Gamma \cap Q(x,r)) \to 0, \quad r \to 0.
\end{equation}

Combining \eqref{eq2.13}, \eqref{eq4.6} and \eqref{eq4.7} and taking into account that $\gamma_{s}$ has the Lusin property (because $\gamma_{s}$ is $1$-Lipschitz) we get
\begin{equation}
\begin{split}
&\varliminf\limits_{r \to 0}\frac{\mu_{0}[\Gamma](Q(x,r))}{r}  \\
&\geq \varliminf\limits_{r \to 0}\fint\limits_{t^{1}_{x}(r)}^{t^{2}_{x}(r)}\chi_{[0,l(\Gamma)]}(\tau)
\operatorname{D}^{0}[\gamma_{s}](\tau)\,d\tau  =
\operatorname{D}^{0}[\gamma_{s}](\gamma^{-1}_{s}(x)) > 0
\quad \hbox{for} \quad \mathcal{H}^{1}-\hbox{a.e.} \quad x \in \Gamma.
\end{split}
\end{equation}

\textit{Step 3.}
Using \eqref{eq264''} we deduce that the measure $\mathcal{H}^{1}\lfloor_{\Gamma}$ is
absolutely continuous with respect to $\mu_{0}[\Gamma]$. Furthermore, $w_{0}[\Gamma](x):=1$ for all $x \in \Gamma$ by 
Theorem \ref{Th2.1} and definition of a $d$-regular sequence of measures given in the introduction. As a result, by \eqref{eq314} we have
\begin{equation}
\label{eq4.9}
w_{k}[\Gamma](x) \geq w_{0}[\Gamma](x) = 1 \quad \hbox{for} \quad \mathcal{H}^{1}-\hbox{a.e.} \quad x \in \Gamma.
\end{equation}
According to Theorem \ref{Th2.1} the sequence of measures $\{\mu_{k}[\Gamma]\}$ is $1$-regular on $\Gamma$.
Hence, we can apply Proposition \ref{Prop2.5} with $d=1$ and then use \eqref{eq264''}, \eqref{eq4.9}. We get the desirable
\begin{equation}
\begin{split}
&\varlimsup\limits_{k \to \infty}\frac{1}{w_{k}[\Gamma](x)}\fint\limits_{Q_{k}(x)}w_{k}[\Gamma](y)\,d\mu_{0}[\Gamma](y)\\
&=\varlimsup\limits_{k \to \infty}\frac{1}{w_{k}[\Gamma](x)}\frac{\mu_{k}[\Gamma](Q_{k}(x))}{\mathcal{H}^{1}\lfloor_{\Gamma}
(Q_{k}(x))}
\frac{\mathcal{H}^{1}\lfloor_{\Gamma}(Q_{k}(x))}{\mu_{0}[\Gamma](Q_{k}(x))}\\
&\le \varlimsup\limits_{k \to \infty}\Bigl(\frac{1}{w_{k}[\Gamma](x)}\frac{C^{1}_{\{\mu_{k}[\Gamma]\}}}{g(x)}\Bigr) \le \frac{C^{1}_{\{\mu_{k}[\Gamma]\}}}{g(x)}
< +\infty \quad \hbox{for} \quad \mathcal{H}^{1}-\hbox{a.e.} \quad x \in \Gamma.
\end{split}
\end{equation}
\end{proof}

Recall definition of a $d$-regular sequence of measures on a closed $d$-thick set $S$ given in the introduction.
Recall also Definition \ref{Def.quasipor} and the notion of $\mathcal{BN}_{\{\mathfrak{m}_{k}\},p,\lambda}$ given in
\eqref{eq1.9}.

\begin{Th}
\label{Th.Lebesguepoint}
Let $n\in \mathbb{N}$, $n \geq 2$, $d \in [n-1,n]$, $p \in (1,\infty)$ and $\lambda \in (0,1)$. Let $S \subset \mathbb{R}^{n}$ be a closed $d$-thick set. Let $\{\mathfrak{m}_{k}\}$
be a $d$-regular sequence of measures on $S$. Assume that the following conditions hold:

{\rm (1)} the set $S$ is $(d,\lambda)$-quasi-porous;

{\rm (2)} for $\mathfrak{m}_{0}$-almost every $x \in S$
\begin{equation}
\label{eq2.65'''}
\varlimsup\limits_{k \to \infty}\frac{1}{w_{k}(x)}\fint\limits_{Q(x,2^{-k})}w_{k}(y)\,d\mathfrak{m}_{0}(y) < +\infty.
\end{equation}
Then, the condition $\mathcal{BN}_{\{\mathfrak{m}_{k}\},p,\lambda}[f]  < +\infty$ implies
\begin{equation}
\label{eq239'}
\lim\limits_{k \to \infty}\fint\limits_{Q(x,2^{-k})}|f(x)-f(y)|\,d\mathfrak{m}_{k}(y)=0 \quad \hbox{for} \quad \mathfrak{m}_{0}-\hbox{a.e.} \quad x
\in S.
\end{equation}
\end{Th}

\begin{proof}

We split the proof into several steps.

\textit{Step 1.}
Using Proposition \ref{Prop23} we get for each $k \in \mathbb{N}_{0}$
\begin{equation}
\label{eq2.64''}
2^{k(p-(n-d))}
\Bigl(\widetilde{\mathcal{E}}_{\mathfrak{m}_{k}}[f](Q_{k}(x))\Bigr)^{p} \le 2^{k(d-n)}\Bigl(2f^{\sharp}_{\{\mathfrak{m}_{k}\}}(x,2^{-k})\Bigr)^{p}
\quad \text{for every} \quad x \in S.
\end{equation}

\textit{Step 2.} By B. Levi theorem  we deduce from \eqref{eq2.64''} and \eqref{eq1.9}
\begin{equation}
\label{eq2.65''}
\begin{split}
&\int\limits_{S}\Bigl[\sum\limits_{k=1}^{\infty}2^{k(p-(n-d))}\chi_{S_{k}(\lambda)}(x)
\Bigl(\widetilde{\mathcal{E}}_{\mathfrak{m}_{k}}[f](Q_{k}(x))\Bigr)^{p}
w_{k}(x)\Bigr]\,d\mathfrak{m}_{0}(x)\\
&=\sum\limits_{k=1}^{\infty}2^{k(p-(n-d))}\int\limits_{S_{k}(\lambda)}\Bigl(\widetilde{\mathcal{E}}_{\mathfrak{m}_{k}}[f](Q_{k}(x))\Bigr)^{p}
w_{k}(x)\,d\mathfrak{m}_{0}(x)
\le 2^{p} \Bigl(\mathcal{BN}_{\{\mathfrak{m}_{k}\},p,\lambda}[f]\Bigr)^{p}.
\end{split}
\end{equation}

\textit{Step 3.} Since the set $S$ is $(d,\lambda)$-quasi-porous and $\mathcal{BN}_{\{\mathfrak{m}_{k}\},p,\lambda}[f] < +\infty$ by
\eqref{eq2.65''} we get
\begin{equation}
\label{eq2.69''}
\sum\limits_{k=1}^{\infty}2^{k(p-(n-d))}w_{k}(x)
\Bigl(\widetilde{\mathcal{E}}_{\mathfrak{m}_{k}}[f](Q_{k}(x))\Bigr)^{p} < +\infty \quad \hbox{for} \quad \mathfrak{m}_{0}-\hbox{a.e.} \quad x \in S.
\end{equation}
In particular, this gives
\begin{equation}
\label{eq2.70''}
\lim\limits_{k \to \infty}2^{k(p-(n-d))}
w_{k}(x)\Bigl(\widetilde{\mathcal{E}}_{\mathfrak{m}_{k}}[f](Q_{k}(x))\Bigr)^{p} = 0 \quad \hbox{for} \quad \mathfrak{m}_{0}-\hbox{a.e.} \quad x \in S.
\end{equation}
Since $p > 1$, $d \in [n-1,n]$ and since $w_{k}(x) \le 2^{k(n-d)}$ for $\mathfrak{m}_{0}$-a.e. $x \in S$ we obtain
\begin{equation}
\label{eq2.71''}
w^{p}_{k}(x) \le 2^{k(p-1)(n-d)}w_{k}(x) \le 2^{k(p-(n-d))}w_{k}(x) \quad \hbox{for} \quad \mathfrak{m}_{0}-\hbox{a.e.} \quad x \in S.
\end{equation}
As a result, combination of \eqref{eq2.70''} and \eqref{eq2.71''} gives
\begin{equation}
\label{eq2.72''}
\lim\limits_{k \to \infty}\Bigl(w_{k}(x)\widetilde{\mathcal{E}}_{\mathfrak{m}_{k}}[f](Q_{k}(x))\Bigr)^{p}=0 \quad \hbox{for} \quad \mathfrak{m}_{0}-\hbox{a.e.} \quad x \in S.
\end{equation}

\textit{Step 4.} Clearly, the condition $\mathcal{BN}_{\{\mathfrak{m}_{k}\}p,\lambda}[f]<+\infty$ implies $f \in L_{1}(\mathfrak{m}_{0})$. Hence, $\mathfrak{m}_{0}$-almost every point $x_{0} \in S$ is a Lebesgue
point of the function $f$ with respect to the measure $\mathfrak{m}_{0}$. This gives
\begin{equation}
\label{eq2.73}
\begin{split}
&\varlimsup\limits_{k \to \infty}\fint\limits_{Q_{k}(x_{0})}|f(x_{0})-f(x)|\,d\mathfrak{m}_{k}(x) \le\\
&\varlimsup\limits_{k \to \infty}\fint\limits_{Q_{k}(x_{0})}|f(x_{0})-f(y)|\,d\mathfrak{m}_{0}(y)+
\varlimsup\limits_{k \to \infty}\fint\limits_{Q_{k}(x_{0})}\fint\limits_{Q_{k}(x_{0})}|f(y)-f(x)|\,d\mathfrak{m}_{0}(y)d\mathfrak{m}_{k}(x)\\
&\le \varlimsup\limits_{k \to \infty}\fint\limits_{Q_{k}(x_{0})}\fint\limits_{Q_{k}(x_{0})}|f(y)-f(x)|\,d\mathfrak{m}_{0}(y)d\mathfrak{m}_{k}(x)
\quad \hbox{for} \quad \mathfrak{m}_{0}-\hbox{a.e.} \quad x \in S.
\end{split}
\end{equation}

\textit{Step 5.}
Since  $w_{k}(x) \geq 1$ for $\mathfrak{m}_{0}$-a.e. $x \in S$ and since $\operatorname{supp}\mathfrak{m}_{k} = S$, $k \in \mathbb{N}_{0}$ we have by \eqref{eq2.65'''} and \eqref{eq2.72''}
\begin{equation}
\label{eq2.75}
\begin{split}
&\varlimsup\limits_{k \to \infty}\fint\limits_{Q_{k}(x_{0})}\fint\limits_{Q_{k}(x_{0})}|f(y)-f(x)|\,d\mathfrak{m}_{k}(x)\,d\mathfrak{m}_{0}(y)\\
&\le \varlimsup\limits_{k \to \infty} \frac{\mathfrak{m}_{k}(Q_{k}(x_{0}))}{w_{k}(x_{0})\mathfrak{m}_{0}(Q_{k}(x_{0}))}\frac{w_{k}(x_{0})}{\mathfrak{m}_{k}(Q_{k}(x_{0}))}\int\limits_{Q_{k}(x_{0})}
\fint\limits_{Q_{k}(x_{0})}|f(y)-f(x)|d\mathfrak{m}_{k}(x)w_{k}(y)\,d\mathfrak{m}_{0}(y)\\
&\le \varlimsup\limits_{k \to \infty}\frac{\mathfrak{m}_{k}(Q_{k}(x_{0}))}{w_{k}(x_{0})\mathfrak{m}_{0}(Q_{k}(x_{0}))}
\varlimsup\limits_{k \to \infty}w_{k}(x_{0})
\widetilde{\mathcal{E}}_{\mathfrak{m}_{k}}[f](Q_{k}(x_{0}))=0, \quad \mathfrak{m}_{0}-\hbox{a.e.} \quad x_{0} \in S.
\end{split}
\end{equation}

Combining \eqref{eq2.73} with \eqref{eq2.75} we get \eqref{eq239'} and complete the proof.

\end{proof}

If $\Gamma \subset \mathbb{R}^{2}$ is a simple rectifiable curve of positive length then for any
$1$-regular on $\Gamma$ sequence of measures $\{\mathfrak{m}_{k}\}_{k \in \mathbb{N}_{0}}$ every measure $\mathfrak{m}_{k}$, $k \in \mathbb{N}_{0}$ is finite on $\Gamma$.
Hence, using \eqref{eq2.17}  we have $f \in L_{1}(\mathfrak{m}_{k_{0}})$ for some $k_{0} \in \mathbb{N}_{0}$
if and only if $f \in L_{1}(\mathfrak{m}_{k})$ for all $k \in \mathbb{N}_{0}$. Furthermore, if
by H\"older inequality if $f \in L_{p}(\mathfrak{m}_{k_{0}})$ for some $k_{0} \in \mathbb{N}_{0}$ and $p \in [1,\infty)$
then $f \in L_{1}(\mathfrak{m}_{k})$ for all $k \in \mathbb{N}_{0}$.
Recall the construction of the extension operator \eqref{eq2.37'} given in the introduction.
Given a Borel function $f \in L_{p}(\mu_{0}[\Gamma])$, we define (we set $k(r):=[\log_{2}r^{-1}]$)
\begin{equation}
\label{eq2.37''''}
F(x)=\operatorname{Ext}_{\Gamma}[f](x):=\sum\limits_{\alpha \in \mathcal{I}}\varphi_{\alpha}(x)\fint\limits_{\widetilde{Q}_{\alpha} \cap \Gamma}f(\widetilde{x})
\,d\mu_{k(r_{\alpha})}[\Gamma](\widetilde{x})
,\quad x \in \mathbb{R}^{2}.
\end{equation}

Now we can formulate the main result of this section.

\begin{Th}
\label{Th2.2}
Let $\Gamma \subset [0,1) \times [0,1)$ be a simple rectifiable curve of positive length. Let $\gamma: [a,b] \to \mathbb{R}^{2}$
be an admissible parametrization of $\Gamma$.
Let $\{\mu_{k}\}=\{\mu_{k}[\Gamma]\}_{k \in \mathbb{N}_{0}}$  be the special $1$-regular sequence of measures
on $\Gamma$.
Given $p \in (1,\infty)$, a function $f: \Gamma \to \mathbb{R}$ belongs to the trace space $W^{1}_{p}(\mathbb{R}^{2})|_{\Gamma}$ if and only if there exists $\lambda_{0}:=\lambda_{0}(\Gamma) \in (0,1]$
such that
\begin{equation}
\begin{split}
&\mathcal{BN}_{\{\mu_{k}\},p,\lambda_{0}}[f]:=\|f|L_{p}(\mu_{0})\|+\Bigl(\sum_{k=1}^{\infty}2^{-k}\int\limits_{\Gamma_{k}(\lambda_{0})}
\Bigl(f^{\sharp}_{\{\mu_{k}\}}(x,2^{-k})\Bigr)^{p}\,d\mu_{k}(x)\Bigr)^{\frac{1}{p}} < \infty.
\end{split}
\end{equation}
Furthermore, for any $\lambda \in (0,\lambda_{0}]$
there  exists a constant $C > 0$ depending only on  $p$, $\lambda$, $C^{1}_{\{\mu_{k}\}}$, $C^{2}_{\{\mu_{k}\}}$ such that
\begin{equation}
C^{-1}\mathcal{BN}_{\{\mu_{k}\},p,\lambda}[f] \le \|f|W_{p}^{1}(\mathbb{R}^{2})|_{\Gamma}\| \le C \mathcal{BN}_{\{\mu_{k}\},p,\lambda}[f]
\end{equation}
and the operator $\operatorname{Ext}_{\Gamma}: W_{p}^{1}(\mathbb{R}^{2})|_{\Gamma} \to W_{p}^{1}(\mathbb{R}^{2})$ defined in
\eqref{eq2.37''''} is linear, bounded and $\operatorname{Tr}|_{\Gamma} \circ \operatorname{Ext}_{\Gamma} = \operatorname{Id}$
on the space $W_{p}^{1}(\mathbb{R}^{2})|_{\Gamma}$.
\end{Th}

\begin{proof}
Since $\Gamma$ is a path-connected set and since $\Gamma$ is compact we get existence of some constant $c > 0$  such that
$\mathcal{H}^{1}(Q(x,r) \cap \Gamma) \geq c r$ for and all $x \in \Gamma$ and $r \in (0,1]$.
By Proposition \ref{Prop.curve} and Lemma \ref{Lm.porous} it follows that there is a number $\widetilde{\lambda}_{0}:=\widetilde{\lambda}_{0}(\Gamma) \in (0,1]$
such that the curve $\Gamma$ is $(1,\lambda)$-quasi-porous for any $\lambda \in (0,\widetilde{\lambda_{0}}]$. Hence, combining
Proposition \ref{Prop2.5}, Theorem \ref{Th2.1}, Lemma \ref{Prop2.6} and Theorem \ref{Th.Lebesguepoint} we deduce that for any $\lambda \in (0,\widetilde{\lambda_{0}}]$ the condition
$\mathcal{BN}_{\{\mu_{k}\},p,\lambda}[f] < +\infty$ implies that
\begin{equation}
\notag
\lim\limits_{k \to \infty}\fint\limits_{Q(x,2^{-k})}|f(x)-f(y)|\,d\mu_{k}(y)=0 \quad \hbox{for} \quad \mu_{0}[\Gamma]-\hbox{a.e.} \quad x \in \Gamma.
\end{equation}
On the other hand, it was mentioned in the proof of Lemma \ref{Prop2.6} that the measure $\mathcal{H}^{1}\lfloor_{\Gamma}$ is absolutely continuous
with respect to $\mu_{0}[\Gamma]$. Hence,
$\mathcal{BN}_{\{\mu_{k}\},p,\lambda}[f] < +\infty$ implies that
\begin{equation}
\notag
\lim\limits_{k \to \infty}\fint\limits_{Q(x,2^{-k})}|f(x)-f(y)|\,d\mu_{k}(y)=0 \quad \hbox{for} \quad \mathcal{H}^{1}-\hbox{a.e.} \quad x \in \Gamma.
\end{equation}
Combining this observation  with Theorem \ref{ThA} and Theorem \ref{Th2.1} we complete the proof.

\end{proof}

\section{Example}

Note that despite the fact that for a given simple rectifiable curve $\Gamma \subset \mathbb{R}^{2}$ with positive length
the measures $\mu_{k}[\Gamma]$, $k \in \mathbb{N}_{0}$ constructed above  
have explicit expressions, still a generic simple rectifiable curve $\Gamma \subset \mathbb{R}^{2}$ can be extremely complicated and
functions $\operatorname{D}^{k}[\gamma]$ can oscillate wildly in general. As a result, given a 
function $f \in W_{p}^{1}(\mathbb{R}^{2})|_{\Gamma}$ computations of
the norm $\mathcal{BN}_{\{\mu_{k}[\Gamma]\},p,\lambda}[f]$ can be problematic in practice.

Below we present an illustrative example,
which on the one hand is quite simple for computations, on the other hand it exhibits the typical
effects when oscillations of a given curve $\Gamma$ affect to the behavior of densities $\operatorname{D}^{k}[\gamma]$ 
(for admissible parameterizations $\gamma$).

We restrict ourselves to the case when $\Gamma$ is a graph of some locally Lipschitz nonnegative function. More precisely, we assume that $\gamma_{1}=\operatorname{id}$ on $[0,1]$ and $\gamma_{2}:[0,1] \to \mathbb{R}_{+}$ is a locally Lipschitz (i.e. $\gamma_{2}|_{[a,b]}$ is Lipschitz for each $[a,b] \subset (0,1)$)  function. We set $\Gamma:=\{(x_{1},x_{2}):x_{1} \in [0,1], x_{2}=\gamma_{2}(x_{1})\}$. Note that in order to make our example
interesting our curve $\Gamma$ should satisfy the following requirements:

{\rm a)} the parametrization $\gamma$ is not (globally) Lipschitzian because otherwise we fall into the scope of \cite{Gal};

{\rm b)} the graph $\Gamma$ fails to satisfy the Ahlfors-David 1-regularity condition
because otherwise we fall into the scope of \cite{Ihn}. Hence, $\gamma_{2}$
should oscillate strongly.

Firstly we define
\begin{equation}
\zeta(t):=
\begin{cases}
&2t, \quad t \in [0,2^{-1}],\\
&2-2t, \quad t \in (2^{-1},1]\\
&0, \quad t \notin [0,1].
\end{cases}
\end{equation}
For each $k \in \mathbb{N}_{0}$ we set
\begin{equation}
\psi_{k}(t):=\sum\limits_{m \in \mathbb{Z}}\zeta(2^{k}(t-\frac{m}{2^{k}})), \quad t \in \mathbb{R}.
\end{equation}
Let $\{c_{k}\}=\{c_{k}\}_{k \in \mathbb{N}} \subset [0,1)$ be a sequence  of nonnegative numbers and let $\{n_{k}\}=\{n_{k}\}_{k \in \mathbb{N}}$ be a sequence of nonnegative integer numbers
such that $n_{k} > k$ for every $k \in \mathbb{N}$. We define
\begin{equation}
\label{eq5.3'}
\begin{cases}
&\gamma_{1}(t) = t, \quad t \in [0,1];\\
&\gamma_{2}(t):=\sum\limits_{k=1}^{\infty}c_{k}\chi_{[2^{-k},2^{-k+1})}(t)\psi_{n_{k}}(t), \quad t \in [0,1].
\end{cases}
\end{equation}
In other words, the graph $\Gamma$ looks like a sequence of triangles. The amount of the congruent triangles on the interval $[2^{-k},2^{-k+1})$
equals $2^{n_{k}-k}$ and the height  of every such triangle equals $c_{k}$.

Since $\gamma_{2}$ is a continuous function, the map $\gamma=(\gamma_{1},\gamma_{2}):[0,1] \to [0,1) \times [0,1)$ gives a parametrization
of the simple planar curve. Clearly, $\gamma_{2}$ is \textit{locally Lipschitz} (i.e. $\gamma_{2}$ is Lipschitz on any closed interval $[a,1]$
with $a \in (0,1]$).
It is easy to see that the sequences $\{c_{k}\}$ and $\{n_{k}\}$ can be chosen in such a way that:

{\rm (1)}
\begin{equation}
\label{eq2.36}
\sum\limits_{k=1}^{\infty}c_{k}2^{n_{k}-k} < \infty;
\end{equation}

{ \rm (2)}
\begin{equation}
\label{eq2.38}
\lim\limits_{k \to \infty}c_{k}2^{n_{k}} = +\infty.
\end{equation}

As a typical example one can take $c_{k}:=k^{-\alpha}2^{k-n_{k}}$, $k \in \mathbb{N}_{0}$ with $\alpha > 1$.

From now we assume that our planar curve $\Gamma$ satisfies \eqref{eq2.36}--\eqref{eq2.38}.
The condition \eqref{eq2.36} is necessary and sufficient for the \textit{rectifiability} of $\Gamma$.
This fact together with local Lipschitz property of $\gamma_{2}$ implies that the parametrization
$\gamma$ \textit{is admissible}.
From \eqref{eq2.38} it follows that
\begin{equation}
\label{eq2.37}
\lim\limits_{j \to \infty}2^{j}\sum\limits_{k=j}^{\infty}c_{k}2^{n_{k}-k} = +\infty.
\end{equation}
Clearly \eqref{eq2.37} leads to the distortion of the second inequality in \eqref{eq1.1} with $d=1$ and $S=\Gamma$
(the first inequality in \eqref{eq1.1} always holds true with $S=\Gamma$ and $d=1$). Hence, our curve is
locally Lipschitz and \textit{fails to satisfy the Ahlfors-David 1-regularity condition.}
Finally, \eqref{eq2.38} implies that $\gamma_{2}$ is \textit{not globally Lipschitz.}
Hence, the exact description of the trace space $W_{p}^{1}(\mathbb{R}^{2})|_{\Gamma}$ can not be obtained
by earlier known methods of \cite{Gal}, \cite{Ihn}.

It follows directly from the construction that for any $k \in \mathbb{N}_{0}$ it holds
\begin{equation}
\label{eq5.7}
\mathcal{W}^{k}_{2}[\gamma](t)|\dot{\gamma}_{1}(t)|=1 \quad \hbox{for} \quad \hbox{a.e.} \quad t \in [0,1].
\end{equation}
By \eqref{eq5.3'} we obtain for every $k \in \mathbb{N}_{0}$
\begin{equation}
\label{eq5.8}
|\dot{\gamma}_{2}(t)|=c_{k}2^{n_{k}+1} \quad \hbox{for} \quad \hbox{a.e.} \quad t \in [2^{-k-1},2^{-k}].
\end{equation}
Using \eqref{eq5.7} and \eqref{eq5.8} it is easy to deduce (with the help of elementary geometrical arguments) for every $k,l \in \mathbb{N}_{0}$
\begin{equation}
\label{eq2.46}
\chi_{[2^{-l},2^{-l+1})}(t)\operatorname{D}^{k}[\gamma](t) \approx
\chi_{[2^{-l},2^{-l+1})}(t)\widetilde{\operatorname{D}}^{k}[\gamma](t),
\end{equation}
where we set
\begin{equation}
\label{eqq510}
\chi_{[2^{-l},2^{-l+1})}(t)\widetilde{\operatorname{D}}^{k}[\gamma](t)
:=\begin{cases}
&1 \quad \hbox{if} \quad  2^{-k} \in  (0,2^{-n_{l}}) \cup  (c_{l},1] \hbox{ or } c_{l} \le 2^{-n_{l}};\\
&2^{k}c_{l} \quad \hbox{if} \quad 2^{-k} \in [2^{-n_{l}},c_{l}].\\
\end{cases}
\end{equation}
The corresponding constants in \eqref{eq2.46} do not depend on $\gamma$, $k$, $l$ and $t$.
We also define 
\begin{equation}
\label{eq511}
\widetilde{\mu}_{k}[\Gamma]=\gamma_{\sharp}(\widetilde{\operatorname{D}}^{k}[\gamma]\mathcal{H}^{1}\lfloor_{[0,1]}), \quad k \in \mathbb{N}_{0}.
\end{equation}

We would like to describe informally the main idea of \eqref{eq2.46}. Firstly note that $2^{-n_{k}} < c_{k}$ for all sufficiently big $k \in \mathbb{N}$ 
because $c_{k}2^{n_{k}} \to +\infty$, $k \to +\infty$. This implies that \eqref{eqq510} is correct. 
Fix a cube $Q=Q(x,r)$ with $x \in \Gamma$ and $r \approx 2^{-k}$. In order to obtain a simplified version of 
Theorem \ref{ThA} for this case we do not use Theorem \ref{Th2.1}. Indeed in order to obtain a trace criterion 
it is sufficient to obtain not precisely the special sequence of measures $\{\mu_{k}[\Gamma]\}$ but rather something comparable with it.
Hence, we would like to guess in some sense how can we choose the measures $\widetilde{\mu}_{k}[\Gamma]$, $k \in \mathbb{N}_{0}$ in order $\widetilde{\mu}_{k}[\Gamma]$
to satisfy items (1)--(4) in the corresponding definition of $1$-regular on $\Gamma$ sequence of measures.
Roughly speaking, our weight function $\widetilde{\operatorname{D}}^{k}[\gamma]\approx 1$ in the case when either the side length $r$ is comparable
with the base of the corresponding triangle (where the center of the cube is located) or when the side length of the cube is comparable with
the height of the corresponding triangle. In the case when the heights of the congruent triangles (located on the corresponding dyadic intervals) a much bigger than its
bases and when the side length of the cube under consideration is somewhere between these two numbers we should focus on the intersections of lines parallel
to the first coordinate axe with our curve $\Gamma$.

It remains to describe porous subsets of $\Gamma$. First of all we fix $\lambda \in (0,1)$ and $k \in \mathbb{N}_{0}$.
For each $l \in \mathbb{N}$ we set
\begin{equation}
\begin{split}
&E^{l}_{k}(\lambda):=\Bigl(\gamma_{2}^{-1}((c_{l}-\frac{1-\lambda}{2^{k}},c_{l}])\cap [\frac{1}{2^{l}},\frac{2}{2^{l}})\Bigr) \bigcup
\Bigl(\gamma_{2}^{-1}([0,\frac{1-\lambda}{2^{k}}))\cap [\frac{1}{2^{l}},\frac{2}{2^{l}})\Bigr)\\
&\bigcup \Bigl(\gamma_{2}^{-1}((c_{l+1}-\frac{1-\lambda}{2^{k}},c_{l}])\cap [\frac{1}{2^{l}},\frac{1}{2^{l}}+\frac{1-\lambda}{2^{k}})\Bigr) \\
&\bigcup \Bigl(\gamma_{2}^{-1}((c_{l-1}-\frac{1-\lambda}{2^{k}},c_{l}])\cap (\frac{2}{2^{l}}-\frac{1-\lambda}{2^{k}},\frac{2}{2^{l}}])\Bigr).
\end{split}
\end{equation}
Now for every $\lambda \in (0,1)$ and $k \in \mathbb{N}_{0}$ we define the set as
\begin{equation}
E_{k}(\lambda):=\bigcup\limits_{l=1}^{\infty}E^{l}_{k}(\lambda).
\end{equation}
It is an elementary verification that for every $\lambda \in (0,1)$ and $k \in \mathbb{N}_{0}$
\begin{equation}
\label{eq2.49}
\Gamma_{k}(\lambda)= \{(x_{1},x_{2}) \in \mathbb{R}^{2}: x_{1} \in E_{k}(\lambda), x_{2}=\gamma_{2}(x_{1})\}.
\end{equation}

Let $\{\varphi_{\alpha}\}_{\alpha \in \mathcal{I}}$ be the same as in \eqref{eq2.37''''}. Given a Borel
function $f \in L_{p}(\widetilde{\mu_{0}}[\Gamma])$, we define  
\begin{equation}
\label{eq515}
F(x)=\widetilde{\operatorname{Ext}}_{\Gamma}[f](x):=\sum\limits_{\alpha \in \mathcal{I}}\varphi_{\alpha}(x)\fint\limits_{\widetilde{Q}_{\alpha} \cap \Gamma}f(\widetilde{y})
\,d\widetilde{\mu}_{k(r_{\alpha})}[\Gamma](y)
,\quad x \in \mathbb{R}^{2}.
\end{equation}

Application of Theorem \ref{Th2.2} together with \eqref{eq2.46}, \eqref{eq511}, \eqref{eq2.49}, \eqref{eq515} and Proposition \ref{Prop23}
gives the following criterion:

\textit{Given $p \in (1,\infty)$, a function $f: \Gamma \to \mathbb{R}$ belongs to the trace space $W_{p}^{1}(\mathbb{R}^{2})|_{\Gamma}$ if and only if there is
$\lambda_{0} \in (0,1]$ such that}
\begin{equation}
\operatorname{BN}_{\Gamma,p,\lambda_{0}}[f]:=\Bigl(\int\limits_{0}^{1}|f(\gamma_{2}(t))|^{p}\,dt\Bigr)^{\frac{1}{p}}
+\Bigl(\sum\limits_{k=1}^{\infty}2^{-k}\sum\limits_{l=1}^{\infty}
\int\limits_{E^{l}_{k}(\lambda_{0})}\Bigl(\widetilde{f}^{\sharp}_{k}(t)\Bigr)^{p}\widetilde{\operatorname{D}}_{k}[\gamma](t)\,dt\Bigr)^{\frac{1}{p}} < \infty,
\end{equation}
\textit{where we set for each $k \in \mathbb{N}_{0}$ (we use the shorthand $U_{k}(t):=\gamma^{-1}(Q_{k}(\gamma(t)))$)}
\begin{equation}
\label{eq5.14}
\widetilde{f}^{\sharp}_{k}(t):=\sup\limits_{0 \le k' \le k}2^{3k'}\int\limits_{U_{k'}(t)}\int\limits_{U_{k'}(t)}
|f \circ \gamma(s)-f\circ\gamma(s')|\widetilde{\operatorname{D}}_{k}[\gamma](s)\widetilde{\operatorname{D}}_{k}[\gamma](s')\,dsds', \quad t \in [0,1).
\end{equation}

\textit{Furthermore, for each $\lambda \in (0,\lambda_{0}]$ there is a constant $C>0$ depending only on  $p$ and $\lambda$ such that
\begin{equation}
C^{-1}\operatorname{BN}_{\Gamma,p,\lambda}[f] \le \|f|W_{p}^{1}(\mathbb{R}^{2})|_{\Gamma}\| \le C\operatorname{BN}_{\Gamma,p,\lambda}[f]
\end{equation}
and the operator $\widetilde{\operatorname{Ext}}_{\Gamma}: W_{p}^{1}(\mathbb{R}^{2})|_{\Gamma} \to W_{p}^{1}(\mathbb{R}^{2})$ defined in
\eqref{eq515} is linear, bounded and $\operatorname{Tr}|_{\Gamma} \circ \widetilde{\operatorname{Ext}}_{\Gamma} = \operatorname{Id}$
on the space $W_{p}^{1}(\mathbb{R}^{2})|_{\Gamma}$.}

\end{document}